\documentclass[preprint,12pt,reversenotenum]{elsarticle}

\usepackage{microtype}
\usepackage{amsmath,amssymb,amsthm,mathtools,mathrsfs}
\usepackage{esint}
\usepackage{enumitem}
\usepackage{aliascnt}
\usepackage{xcolor}
\usepackage{hyperref}

\mathtoolsset{showonlyrefs=true}

\hypersetup{
  colorlinks=true,
  linkcolor=blue!45!black,
  citecolor=blue!45!black,
  urlcolor=blue!45!black,
  pdftitle={Quantitative Homogenization and Global Lipschitz Estimates for Elasticity with Periodic Rigid Inclusions},
  pdfauthor={Yuanhang Liu and Yuxin Du},
  pdfkeywords={linear elasticity, rigid inclusions, periodic homogenization, pure traction, Lipschitz estimate}
}
\journal{Journal of Differential Equations}
\biboptions{numbers,sort&compress}
\numberwithin{equation}{section}

\newtheorem{theorem}{Theorem}[section]
\newaliascnt{proposition}{theorem}
\newtheorem{proposition}[proposition]{Proposition}
\aliascntresetthe{proposition}
\newaliascnt{lemma}{theorem}
\newtheorem{lemma}[lemma]{Lemma}
\aliascntresetthe{lemma}
\newaliascnt{corollary}{theorem}

\aliascntresetthe{corollary}
\theoremstyle{definition}
\newaliascnt{definition}{theorem}

\aliascntresetthe{definition}
\theoremstyle{remark}
\newaliascnt{remark}{theorem}
\newtheorem{remark}[remark]{Remark}
\aliascntresetthe{remark}

\usepackage[nameinlink,noabbrev]{cleveref}

\crefname{theorem}{Theorem}{Theorems}
\Crefname{theorem}{Theorem}{Theorems}
\crefname{proposition}{Proposition}{Propositions}
\Crefname{proposition}{Proposition}{Propositions}
\crefname{lemma}{Lemma}{Lemmas}
\Crefname{lemma}{Lemma}{Lemmas}
\crefname{corollary}{Corollary}{Corollaries}
\Crefname{corollary}{Corollary}{Corollaries}
\crefname{definition}{Definition}{Definitions}
\Crefname{definition}{Definition}{Definitions}
\crefname{remark}{Remark}{Remarks}
\Crefname{remark}{Remark}{Remarks}
\crefname{section}{Section}{Sections}
\Crefname{section}{Section}{Sections}
\crefname{equation}{equation}{equations}
\Crefname{equation}{Equation}{Equations}

\newcommand{\R}{\mathbb R}
\newcommand{\Z}{\mathbb Z}
\newcommand{\T}{\mathbb T}
\newcommand{\sym}{\operatorname{sym}}
\newcommand{\diver}{\operatorname{div}}
\newcommand{\dist}{\operatorname{dist}}

\newcommand{\supp}{\operatorname{supp}}

\newcommand{\Id}{\mathrm I}
\newcommand{\Ahom}{A^{\mathrm{hom}}}
\newcommand{\avgint}{\fint}
\newcommand{\eps}{\varepsilon}
\newcommand{\ee}{e}
\newcommand{\calR}{\mathcal R}

\newcommand{\norm}[1]{\left\lVert #1\right\rVert}
\newcommand{\abs}[1]{\left\lvert #1\right\rvert}
\newcommand{\set}[1]{\left\{#1\right\}}
\newcommand{\pair}[2]{\left\langle #1,#2\right\rangle}
\begin{document}

\begin{frontmatter}

\title{Quantitative Homogenization and Global Lipschitz Estimates
for Elasticity with Periodic Rigid Inclusions}
\author[aff1,aff2]{Yuanhang Liu}
\ead{liuyh23@mails.tsinghua.edu.cn}
\author[aff1]{Yuxin Du\corref{cor1}}
\ead{duyx23@mails.tsinghua.edu.cn}

\cortext[cor1]{Corresponding author}

\address[aff1]{Qiuzhen College, Tsinghua University,
Beijing 100084, China}

\address[aff2]{Institute for AI Industry Research (AIR),
Tsinghua University, Beijing, China}

\begin{abstract}
We consider the pure traction problem for linear elasticity in a bounded
domain containing an $\varepsilon$-periodic array of perfectly rigid
inclusions. For compatible $L^2$ boundary tractions on a Lipschitz domain,
we prove an $O(\varepsilon^{1/2})$ error estimate in $H^1$ for a smoothed
first-order two-scale expansion. If the exterior and inclusion boundaries
are $C^{2,\beta}$, we also obtain a global Lipschitz estimate, uniform in
$\varepsilon$, for compatible $C^{0,\sigma}$ tractions, where
$0<\sigma<\beta$. The displacement is normalized modulo rigid motions,
and the estimate controls the full gradient in the matrix and the
inclusions. The boundary argument uses a local Neumann comparison estimate
whose constant is independent of the choice of inclusions retained in
cells cut by an artificial localization boundary. We prove this estimate
by combining cellwise rigidification, divergence-free stress extensions,
and duality, and then apply a boundary Campanato iteration and local
estimates at the microscopic scale.
\end{abstract}


\end{frontmatter}
\section{Introduction}

\subsection{Geometric setting and the problem set-up}

We begin by fixing the geometry and the variational formulation of the microscopic problem.

Let $d\ge2$, let $\Omega\subset\R^d$ be a bounded connected Lipschitz domain, and let $ Y=[-1/2, 1/2 )^d $ be the unit cell. Let $D$ be a $C^{1,\alpha}$ domain with connected boundary compactly contained in $Y$ for some $\alpha \in (0,1)$, and its distance to $ \partial Y$ is $ \dist(D,\partial Y)> \delta_0  $.

 For any $\varepsilon \in (0,1)$, define
\begin{equation}\label{jepsi}
    J_{\varepsilon}:=\{ \mathbf{m} \in \mathbb{Z}^d : \varepsilon (Y+\mathbf{m}) \subset \Omega \},
\end{equation}
and the inclusions set
\begin{equation}\label{def0D}
    D_{\varepsilon}:= \bigcup_{\mathbf{m}\in J_{\varepsilon}} \varepsilon (D+\mathbf{m}) .
\end{equation}
Assume the Lam\'e constants satisfy $\mu>0$ and $d\lambda+2\mu>0$, the Lam\'e system is
\begin{equation}
    \mathcal{L}_{\lambda,\mu} u := \mu \Delta u + (\lambda +\mu) \nabla \diver u.
\end{equation}
Denote the elliptic tensor
\begin{equation}
    \mathbb C E=\lambda\operatorname{tr}(E)\Id+2\mu E, \quad  \forall E \in\R^{d\times d}_{\sym},
\end{equation}
and let $ \ee $ be the symmetrized differential operator:
\begin{equation}
    \ee(v)=\frac12(\nabla v+\nabla v^T)= \frac{1}{ 2} (\partial_i v_j + \partial_j v_i ),
\end{equation}
here the superscript $ T $ is the transpose of matrices. Then the Lam\'e system can also be written as
\begin{equation}
    \mathcal{L}_{\lambda,\mu} u = \diver \mathbb{C} \ee(u).
\end{equation}
It is clear that $ \min\{2\mu,2\mu+d\lambda\}|E|^2 \le \mathbb{C}E:E \le (2\mu+d|\lambda|)|E|^2 $, for all $ E\in\R^{d\times d}_{\sym} $. The corresponding boundary traction at $ \partial\Omega $ is
\begin{equation}
    \frac{\partial u}{\partial \nu_{(\lambda,\mu)} }\bigg|_{\partial \Omega } : = \mathbb{C} \ee(u) N = \lambda (\diver u ) N + 2\mu \ee(u) N,
\end{equation}
here $ N $ is the outward unit normal vector to $ \partial\Omega $.

Let
\begin{equation}
    \calR=\set{r(x)=a+Bx:a\in\R^d,\ B^T=-B}
\end{equation}
be the space of rigid motions. It is clear that $ \mathcal{R} $ is a linear space with dimension $ d(d+1)/2 $. Denote $H^{-1/2}_{\calR}(\partial\Omega)$ as the subspace of $ H^{-1/2}(\partial\Omega) $ that is orthogonal to $ \mathcal{R} $:
\begin{equation}
     H^{-1/2}_{\calR}(\partial\Omega)
 =\set{g\in H^{-1/2}(\partial\Omega;\R^d):
        \pair g r_{\partial\Omega}=0\ \text{for every }r\in\calR},
\end{equation}
and denote also $L^2_{\calR}(\partial\Omega)$ as the subspace of $ L^2(\partial\Omega) $ that is orthogonal to $ \mathcal{R} $:
\begin{equation}
     L^2_{\calR}(\partial\Omega)
 =\set{g\in L^2(\partial\Omega;\R^d):
        \int_{\partial\Omega}g\cdot r\,dS=0\ \text{for every }r\in\calR}.
\end{equation}
Consider the following problem
\begin{equation}\label{Lame_eq_main}
    \left\{ \begin{array}{ll}
        \mathcal{L}_{\lambda,\mu} u_\eps = 0, & \mathrm{in}\ \Omega \setminus D_\eps,  \\
        \ee(u_\eps) = 0,  & \mathrm{in}\ D_\eps,\\
        \frac{\partial u}{\partial \nu_{(\lambda,\mu)} }\big|_{\partial\Omega} = g \in H^{-1/2}_{\calR}(\partial\Omega),
    \end{array}
    \right.
\end{equation}
in the sense that: let
\begin{equation}
    V_\eps:=\set{v\in H^1(\Omega;\R^d):
        \ee(v)=0\ \text{in }D_\eps,\quad
        \int_{\partial\Omega}v\cdot r=0\ \text{for every }r\in\calR}.
\end{equation}
then for each $ \varphi\in V_\eps $,
\begin{equation}\label{eq:micro-weak-intro}
  \int_{\Omega \setminus D_\eps }\mathbb C\ee(u_\eps):\ee(\varphi)
  =\pair g\varphi_{\partial\Omega}
  \qquad(\varphi\in V_\eps).
\end{equation}
The existence and uniqueness of the solution in $ V_\eps $ of \eqref{Lame_eq_main} follows from Korn's inequality. Denote 
\begin{equation}
    \widetilde V_\eps : = \set{v\in H^1(\Omega;\R^d):
        \ee(v)=0\ \text{in }D_\eps },
\end{equation}
it then follows from $ g\in H_\mathcal{R}^{-1/2} (\partial \Omega) $ that \eqref{eq:micro-weak-intro} still holds for each $ \varphi\in \widetilde V_\eps $




\subsection{Motivation and Background}


Periodic composites with highly stiff inclusions provide a basic model for the interaction between macroscopic elastic deformation and microscopic rigid-body motion. In the infinite-contrast limit, the inclusions become perfectly rigid: the displacement remains defined throughout the body, but on each inclusion it is constrained by
\[
\ee(u_\eps)=0,
\]
and therefore is an a priori unknown rigid motion. This differs essentially from perforated or clamped-inclusion problems, where the displacement is either absent or prescribed inside the inclusions.

The qualitative homogenization of periodic media and perforated domains is by now classical; see, for instance, \cite{MR578345,MR652509,cioranescu1999introduction,jikov1994homogenization}, and two-scale convergence provides a general framework for identifying homogenized limits \cite{doi:10.1137/0523084}. On the quantitative side, the compactness method of Avellaneda and Lin \cite{avellaneda1987compactness} initiated a systematic theory of uniform regularity for periodic elliptic equations and systems. Boundary estimates, including the Neumann theory relevant here, were subsequently developed in \cite{kenig2013neumann,shen2017boundary,shen2018periodic}. For Lam\'{e} systems in perforated media with clamped holes, related homogenization and boundary estimates were obtained in \cite{jing2021layerPotentials}.

Rigid inclusions may also be viewed as the hard-inclusion endpoint of high-contrast elasticity. In the periodic setting, Shen \cite{shen2021largeScale} established large-scale interior Lipschitz estimates uniformly with respect to the contrast, including the stiff limit. Fu and Jing established quantitative convergence to the soft- and hard-inclusion limits for periodic elastostatic systems \cite{fuJing2024uniform}; more recently, they obtained contrast-uniform convergence rates and global Lipschitz estimates for bounded-domain Dirichlet problems involving high-contrast elliptic systems \cite{fuJing2025rates}. For other constrained models, \cite{signing2024rigid} studied the mixed clamped--traction problem with a Robin-type interaction. 

The work most closely related to the present paper is that of Duerinckx and Gloria \cite{duerinckxGloria2022suspensions}, who developed quantitative stochastic homogenization and large-scale regularity for random suspensions of rigid particles in steady Stokes flow. As they point out, their methods also apply to compressible and incompressible elasticity with rigid inclusions, including the periodic setting. Accordingly, neither the exact rigidity constraint nor stress extension, local surgery, flux-corrector constructions, and large-scale interior regularity are claimed here as new by themselves.

The distinction lies in the boundary problem and in the endpoint regularity estimate. The regularity theory in \cite{duerinckxGloria2022suspensions} is an interior large-scale theory, formulated through spatial averages of the full gradient on scales above a minimal radius, while their bounded-domain quantitative result is posed with homogeneous Dirichlet data. It does not address regularity up to an exterior traction boundary. In the present paper, we consider a bounded elastic body under pure traction and establish the deterministic pointwise estimate
\begin{equation}
    \|\nabla u_\varepsilon\|_{L^\infty(\Omega)}
\leq C\|g \|_{C^{0,\sigma}(\partial\Omega)},
\end{equation}
uniformly in $ \varepsilon$, throughout both the elastic matrix and the rigid inclusions and up to the exterior boundary. Thus the estimate is not only large-scale: it remains valid down to and below the particle scale.




\subsection{Main results}

We state the main results of the paper in this section.

One of the main tasks of this article is to prove the homogenized problem of \eqref{Lame_eq_main} is 
\begin{equation}\label{lame_eq_homo}
    \left\{ \begin{array}{ll}
        \diver \Ahom \ee(u_0) = 0, & \mathrm{in}\ \Omega ,  \\
        \Ahom\ee(u_0)N |_{\partial\Omega} = g \in H^{-1/2}_{\calR}(\partial\Omega),
    \end{array}
    \right.
\end{equation}
here the effective tensor is defined by
\begin{equation}\label{eq:Ahom-intro}
  (\Ahom)_{ij}^{\alpha\beta }
  =\int_{Y}\mathbb C\ee(\chi^\alpha_i + P_i^\alpha) : \ee(\chi^\beta_j + P_j^\beta),
\end{equation}
here $ \chi_i^\alpha  $ is the solution of the cell problem below \eqref{cell_pb}, and $ P^\alpha_i (y) = y_i e^\alpha $, $ e^\alpha = (0,\cdots,1, \cdots,0) $ with $1$ in the $ \alpha $ position. It is clear that $ (\Ahom)_{ij}^{\alpha\beta } = (\Ahom)_{i\beta}^{\alpha j }  = (\Ahom)_{\alpha j}^{i\beta } $ and $ \Ahom $ is symmetric, i.e., $ (\Ahom)_{ij}^{\alpha\beta } = (\Ahom)_{j i}^{\beta \alpha} $.

\begin{theorem}\label{thm:main-hom}
The homogenized equation \eqref{lame_eq_homo} has a unique solution in
\begin{equation}
     H^1_{\calR}(\Omega;\R^d)
 =\set{v\in H^1(\Omega;\R^d):
          \int_{\partial\Omega}v\cdot r=0\ \text{for every }r\in\calR}.
\end{equation}
Let $g\in H^{-1/2}_{\calR}(\partial\Omega)$ and let $u_\eps\in V_\eps$ be the unique solution of \eqref{Lame_eq_main}, let $ u_0\in H^1_{\calR}(\Omega) $ be the unique solution of \eqref{lame_eq_homo}, then
\begin{equation}
    u_\eps \rightarrow u_0 \ \mathrm{in}\ L^2(\Omega;\R^d), \quad \mathrm{and}\ u_\eps \rightharpoonup u_0 \ \mathrm{in}\ H^1(\Omega).
\end{equation}
\end{theorem}



\begin{theorem}\label{thm:main-rate}
There exists a constant $ C = C(\lambda, \mu,d,\Omega, D) $, such that for $g\in L^2_{\calR}(\partial\Omega)$ and $ u_\eps \in V_\eps, u_0 \in  H^1_{\calR}(\Omega)$ solutions of \eqref{Lame_eq_main} and \eqref{lame_eq_homo} respectively, we have
\begin{equation}\label{444}
    \begin{aligned}
        \norm{u_\eps-u_0 - \eps \chi_i^\alpha(\frac{\cdot}{\eps}) S_\eps \big(\eta_\eps \frac{\partial u_0^\alpha}{\partial x_i} \big)}_{H^1(\Omega)}
  &\le C\eps^{1/2}\norm g_{L^2(\partial\Omega)},
    \end{aligned}
\end{equation}
here $ S_\eps $ is the smoothing operator defined as \eqref{epsi_smoother}, $ \eta_\eps $ is a cutoff function defined as \eqref{cut-off_def}.
\end{theorem}

\begin{remark}
 The cutoff removes the oscillatory corrector near the exterior boundary. The resulting expansion is not generally rigid in the inclusions. In the proof, a cellwise rigidification produces an admissible comparison field with the same exterior trace, and its difference from the expansion is estimated separately. The boundary-layer estimate for the homogenized Neumann problem then gives \eqref{444}. In particular, Theorem \ref{thm:main-rate} does not require a global $H^2$ estimate for $u_0$.
\end{remark}



\begin{theorem}\label{thm:main-global}
Assume in addition $\partial\Omega,\partial D\in C^{2,\beta}$, let $0<\sigma<\beta$, and let
\[
 g\in C^{0,\sigma}(\partial\Omega;\R^d)
       \cap H^{-1/2}_{\calR}(\partial\Omega).
\]
Then the solution of \eqref{Lame_eq_main} satisfies
\begin{equation}\label{eq:main-global}
  \norm{\nabla u_\eps}_{L^\infty(\Omega)}
  \le C\norm g_{C^{0,\sigma}(\partial\Omega)}
  \qquad(0<\eps\le1).
\end{equation}
\end{theorem}

\begin{remark}
On each inclusion, $u_\varepsilon(x)=a_{\varepsilon,k}+B_{\varepsilon,k}x$ with $B_{\varepsilon,k}^T=-B_{\varepsilon,k}$. Thus \eqref{eq:main-global} also bounds every particle rotation. Korn's inequality controls the global $L^2$ norm of $\nabla u_\varepsilon$ after normalization, but does not by itself give a bound for the individual matrices $B_{\varepsilon,k}$ that is uniform in $\varepsilon$. The normalization is essential, since adding a global rigid rotation preserves the equation and the traction while changing the full gradient.

\end{remark}


\subsection{Main contribution}

The principal regularity contribution of this paper is a global, full-gradient Lipschitz estimate compatible with both exact rigidity and pure traction. Two additional mechanisms are needed to reach this endpoint.

The first concerns localization at the exterior boundary. A boundary Campanato argument requires local Neumann comparison problems on domains whose artificial boundaries may intersect periodic cells. Pure traction preserves the rigid-motion kernel, so the localized boundary datum must satisfy
\begin{equation}
    \langle h,r\rangle_{\partial G}=0
\qquad
\text{for every }r\in\mathcal R,
\end{equation}
equivalently, both the total force and the total torque must vanish. Arbitrarily restricting the microscopic traction to the localization domain does not make this compatibility transparent.

We resolve this difficulty by extending the microscopic stress through the rigid particles to a symmetric divergence-free tensor $\widehat\sigma_\varepsilon$ on the whole localization domain and prescribing
\begin{equation}
    h_G=\widehat\sigma_\varepsilon N
\qquad\text{on }\partial G.
\end{equation}
For every rigid motion $ r\in \mathcal{R} $,
\begin{equation}
    \langle h_G,r\rangle_{\partial G}
=\int_{\partial G}\widehat\sigma_\varepsilon N \cdot r  =0,
\end{equation}
Thus the local traction automatically satisfies the required force and torque compatibility conditions. Combined with cellwise rigidification, this construction yields a local homogenized Neumann comparison estimate whose constants are uniform with respect to the retention or deletion of cells intersected by the artificial boundary. The novelty here is not the stress extension alone, but its use as an equilibrated localization device for the pure-traction boundary problem.

The second mechanism concerns the passage from large-scale control to a pointwise estimate for the full gradient. On every connected rigid inclusion,
\begin{equation}
    u_\varepsilon(x)
=a_{\varepsilon,k}
+B_{\varepsilon,k}x,
\qquad
B_{\varepsilon,k}^{T}=-B_{\varepsilon,k}.
\end{equation}
The elastic energy only controls $\ee(u_\varepsilon)$ and is therefore insensitive to the matrices $B_{\varepsilon,k}$. Large-scale Campanato control does not by itself provide a point-wise bound for these microscopic rotations.

At the microscopic scale, we combine local interface and boundary regularity, uniform $W^{1,\infty} $ bounds for the periodic cell correctors, and a non-degeneracy estimate for rigid affine maps. These ingredients give a scale-invariant $L^2$-to-$L^\infty$ estimate on patches of size comparable with \(\varepsilon\) and control every $B_{\varepsilon,k}$ individually. The Campanato iteration first controls the displacement oscillation and the corrected-affine slopes down to $r\simeq\varepsilon$; the microscopic estimate then converts this information into the uniform full-gradient bound
\begin{equation}
    \|\nabla u_\varepsilon\|_{L^\infty(\Omega)}
\leq
C\|g\|_{C^{0,\sigma}(\partial\Omega)}.
\end{equation}
The resulting estimate controls simultaneously the deformation in the elastic matrix, the unknown rotations of all rigid inclusions, and the behavior up to the exterior pure-traction boundary.

\subsection{Strategy of the proof}

We first identify the homogenized operator. Since the physical stress $ \mathbb C\ee(u_\varepsilon)$ is naturally defined only in the elastic matrix, we extend it through every rigid inclusion to obtain a symmetric tensor field $\widehat\sigma_\varepsilon$ on the whole domain such that
\begin{equation}
    \diver\widehat\sigma_\varepsilon=0
\quad\text{in }\Omega,
\qquad
\widehat\sigma_\varepsilon N=g
\quad\text{on }\partial\Omega.
\end{equation}
The corresponding cell stresses are extended in the same way. These divergence-free oscillating stresses provide a compensated identification of the weak limit and yield
\begin{equation}
    \sigma_0=\Ahom \ee(u_0),
\end{equation}
together with the homogenized pure-traction problem.

To establish the first-order quantitative approximation, we introduce cellwise rigidification operators that modify a comparison field only near an inclusion, make it rigid inside the particle, and leave it unchanged near the boundary of the cell, and extend the cell stress through the inclusion as a symmetric divergence-free tensor with the correct resultant. After subtracting its average, the extended stress admits a standard skew-symmetric flux corrector. Together with smoothing, an exterior cutoff, and the $L^2$-Neumann boundary-layer estimate for the homogenized solution, these constructions yield the first-order $H^1$ error of order $\mathcal O(\varepsilon^{1/2})$.

For the interior Lipschitz estimate, we use a corrected-affine Campanato scheme. The appropriate comparison fields are
\begin{equation}
    b+M(x-x_0)
+\varepsilon\chi_i^\alpha(x/\varepsilon)M_i^\alpha.
\end{equation}
An adapted Caccioppoli inequality and a local $ \Ahom $-harmonic approximation give a one-step decay of the corresponding displacement excess. A coercivity property of the corrected-affine parametrization controls the matrices selected at successive scales. Iteration then bounds both the excess and the affine slopes down to a scale comparable with $\eps$.

For the boundary estimate, we use the equilibrated localization mechanism described above to construct a local homogenized Neumann comparison problem, uniformly with respect to the cells intersecting the artificial boundary. The boundary excess measures both the deviation of \(u_\varepsilon\) from an affine displacement and the mismatch between the prescribed traction and the corresponding homogenized affine conormal. Constant-coefficient Neumann regularity gives decay of this excess, while the local homogenized comparison transfers the decay to the microscopic solution with an error of order $(\varepsilon/r)^{1/2}$. A discrete Campanato iteration then controls the affine slopes and the displacement oscillation down to $r\simeq\varepsilon$.

Finally, the microscopic \(L^2\)-to-\(L^\infty\) estimate described above converts the scale-\(\varepsilon\) Campanato control into a pointwise bound for the full gradient. This closes the argument in both phases and up to the exterior boundary.

\subsection{Organization of the article}

The paper is organized as follows. Section~2 identifies the homogenized tensor and proves qualitative homogenization, then develops the compatible reconstruction, whole-cell stress extension, flux corrector, and the first-order \(H^1\) convergence rate. Section~3 establishes the Caccioppoli inequalities and the local \(\Ahom\)-harmonic comparison estimates needed to follow the approximation scheme. Section~4 proves the large-scale interior and boundary estimates and completes the microscopic closure needed for the global Lipschitz bound.

\section{Qualitative and quantitative homogenization}

\subsection{The cell problem and modification of the stress}

We first derive the cell problem:

Adopt the following two-scale ansatz:
\begin{equation}\label{ansatz_inside}
     u_\varepsilon(x) = u_0(x,x/\varepsilon) + \varepsilon u_1(x,x/\varepsilon)  + \cdots  \quad\mathrm{in} \ \Omega,
\end{equation}
where, for each $j\in \mathbb{N}$, the function $u_j$ is $Y$-periodic in the fast variable $y:=x/\varepsilon$ and depends on the slow variable $x$. For $ j\ge 1 $, we add the normalization $ \int_Y u_j = 0 $.

Substituting \eqref{ansatz_inside} into \eqref{Lame_eq_main}, using the chain rule
\begin{equation}
    \nabla = \nabla_x +\varepsilon^{-1 }\nabla_y ,
\end{equation}
and equating coefficients of like powers of $\varepsilon$, yield
\begin{equation}\label{240}
    \begin{aligned}
        &\diver_y \mathbb{C} \ee_y(u_0) = 0, \\
        & \diver_y \mathbb{C} (\ee_y(u_1) + \ee_x(u_0))+ \diver_x\mathbb{C} \ee_y(u_0)=0,\\
    \end{aligned}
    \qquad \mathrm{in}\ \Omega\times Y\setminus D,
\end{equation}
and
\begin{equation}\label{248}
    \begin{aligned}
        & \ee_y(u_0) = 0, \\
        & \ee_x(u_0)+ \ee_y(u_1) = 0,
    \end{aligned}\qquad \mathrm{in}\ \Omega\times D,
\end{equation}
with the transmission conditions on $\partial D$:
\begin{equation}\label{255}
    \begin{aligned}
        &\int_{\partial D} \mathbb{C} \ee_y(u_0) N \cdot r= 0, \\
        &  \int_{\partial D} \mathbb{C} (\ee_y(u_1)+\ee_x(u_0)) N \cdot r= 0,
    \end{aligned}\qquad \forall r\in \calR.
\end{equation}
By the first equation of \eqref{248} and \eqref{255}, test \eqref{240} with $ u_0 $ on $ Y\setminus D $ gives
\begin{equation}
    \int_{Y\setminus D} \mathbb{C} \ee_y(u_0): e_y(u_0) = 0,
\end{equation}
hence $ e_y(u_0) \equiv 0 $ on $ \Omega\times Y $, combine with the fact that $ u_0 $ is $ Y $-periodic, we obtain $ u_0(x,y) = u_0(x) $. The second equation of \eqref{240}, \eqref{248} and \eqref{255} then imply that:
\begin{equation}
    u_1(x, y) = \chi_{i}^\alpha(y) (e^\alpha\otimes e^i:\ee(u_0)),
\end{equation}
here $ \chi_{i}^\alpha\in H^1(\T^d; \R^d) $, is the solution of
\begin{equation}\label{cell_pb}
    \left\{ \begin{array}{ll}
        \diver \mathbb{C} \ee(\chi_{i}^{\alpha}+ P_i^\alpha ) = 0, & \mathrm{in} \ Y\setminus D, \\
        \ee(\chi_{i}^\alpha+  P_i^\alpha) =0, & \mathrm{in} \ D, \\
         \int_{\partial D} \mathbb{C} \ee(\chi_{i}^\alpha +P_i^\alpha)N \cdot r = 0, & \forall r\in \calR, \\
         \int_{Y} \chi_{i}^\alpha =0 ,\\
    \end{array}
    \right.
\end{equation}
here $ P^\alpha_i (y) = y_i e^\alpha $, $ e^\alpha = (0,\cdots,1, \cdots,0) $ with $1$ in the $ \alpha $ position. It can be checked that $ \chi_i^\alpha = \chi_\alpha^i $ since $ \ee(P_i^\alpha ) = \ee(P_\alpha^i) $, hence $ u_1 $ can be written in the form of
\begin{equation}
    u_1(x, y) = \chi_{i}^\alpha(y) \frac{\partial u_0^\alpha}{\partial x_i} (x),
\end{equation}


We now state the well-posedness and energy estimate of the cell problem \eqref{cell_pb}, the proof of which is clear.

\begin{proposition}\label{prop:cell-corrector}
For each $1\le j,\alpha\le d $, there exists a unique $ \chi_j^\alpha\in H^1(\T^d; \R^d)$, solving \eqref{cell_pb}, and we have the following energy estimate:
\begin{equation}\label{eq:cell-bound}
  \norm{\chi_j^\alpha}_{H^1(Y)} +\norm{ \ee(\chi_j^\alpha+P_j^\alpha)}_{L^2(Y \setminus D)} \le C .
\end{equation}
\end{proposition}

The existence and uniqueness of solution in $ H_\mathcal{R}^1(\Omega;\R^d) $ of \eqref{lame_eq_homo} follows immediately from Korn's inequality and the following:

\begin{proposition}[Effective tensor]\label{prop:Ahom}
The homogenized elasticity tensor $\Ahom$ is uniformly positive definite. More precisely, there exists a constant $C >1$, such that
\begin{equation}\label{eq:Ahom-ellipticity}
  C^{-1}\abs E^2\le\Ahom E:E\le C \abs E^2,
  \qquad \forall E\in\R^{d\times d}_{\sym}.
\end{equation}
\end{proposition}

\begin{proof}

For each $ E\in\R^{d\times d}_{\sym} $,
\begin{equation}
    \begin{aligned}
        \Ahom E:E & = \int_Y \mathbb{C} \ee(\chi^\alpha_i + P_i^\alpha)E_i^\alpha :\ee(\chi^\beta_j + P_j^\beta)E_j^\beta  \\
        & \ge C^{-1} \int_Y \big| E+\ee(\chi^\alpha_i)E_i^\alpha \big|^2 \ge C^{-1} |E|^2,
    \end{aligned}
\end{equation}
since $ \int \ee(\chi^\alpha_i)E_i^\alpha:E = 0  $. The upper bound follows from \eqref{eq:cell-bound}.
\end{proof}



We next introduce the cell-wise rigidification operators that convert a general comparison field into an admissible one without changing it near the boundary of the cell.

\begin{lemma}\label{lem:rigidification}
There exists a bounded linear map $\mathcal{P}:H^1(Y;\R^d)\longrightarrow H^1(Y;\R^d)$, such that for each $ v\in H^1(Y;\R^d) $,
\begin{equation}
    \ee(\mathcal{P}v ) = 0 \ \mathrm{in}\ D, \quad \mathrm{and} \  \mathcal{P}(v) = v \ \text{in a neighborhood of }\partial Y,
\end{equation}
and
\begin{equation}
    \norm{v- \mathcal{P}v}_{H^1(Y)} \le C\norm{\ee(v)}_{L^2(D)}.
\end{equation}
\end{lemma}

\begin{proof}
    For each $ v\in H^1(Y;\R^d) $, let $ r_v $ be the $L^2(D)$-orthogonal projection onto $\calR$, then Korn's inequality gives
    \begin{equation}
        \|v-r_v\|_{H^1(D) } \le C \norm{\ee(v)}_{L^2(D)}.
    \end{equation}
    Let $ E : H^1(D)\rightarrow H^1(Y) $ be a bounded extension operator such that for each $ u\in H^1(D) $, $ Eu|_D = u $, and $ \|Eu\|_{H^1(Y)} \le C \|u\|_{H^1(D)} $. Take cutoff function $ \zeta \in C_c^\infty (Y;[0,1] ) $ such that $ \zeta\equiv 1 $ on a neighborhood of $ D $, define:
    \begin{equation}
        \mathcal{P} (v):= v -\zeta E(v- r_v),
    \end{equation}
the properties of $ \mathcal{P} $ are easy to verify.
\end{proof}

\begin{lemma}
    \label{lem:rigid-extension}
There exists a bounded linear map $ \mathcal{Q}: H^1(D;\R^d)\longrightarrow H^1(Y;\R^d) $, such that there exist $\zeta\in C_c^\infty(Y)$ and a connected subset $ \widetilde D\Subset Y $, for each $ v\in H^1(D;\R^d) $,
\begin{equation}
    \mathcal{Q} v =v\ \mathrm{in}\ D, \quad \supp(\mathcal{Q}v )\subset \widetilde D , \quad \|\mathcal{Q} v \|_{H^1(Y;\R^d)} \le C \| v \|_{H^1(D;\R^d)} ,
\end{equation}
and for each $ r\in\calR $,
\begin{equation}
    \mathcal{Q} (r) =  \zeta r,
\end{equation}
\end{lemma}

\begin{proof}
    For each $ v\in H^1(D;\R^d) $, let $ r_v $ be the $L^2(D)$-orthogonal projection onto $\calR$, let
    \begin{equation}
        \mathcal{Q}(v) := \zeta r_v + \zeta E(v-r_v),
    \end{equation}
    with $  E $ being the same extension operator and $ \zeta $ be the same cutoff function in proof of Lemma \ref{lem:rigidification}. The properties of $ \mathcal{Q} $ are easy to verify.
\end{proof}

The cell stress $ \mathbb{C} (E+\ee(\chi_E)) $ is initially only divergence free in the elastic matrix \(Y\setminus D\). This is insufficient for the construction of the flux corrector. We therefore extend the cell stress through \(D\) as a symmetric divergence-free tensor whose action is compatible with the rigid-motion kernel. The average of the extended stress is precisely the homogenized stress, and its mean-zero part admits a periodic skew-symmetric flux corrector.

\begin{proposition}[Divergence-free whole-cell stress]\label{prop:whole-cell-stress}
For each $1\le j, \alpha\le d$ there exists
\[
  \widehat\Sigma_j^\alpha\in L^2_\#(Y;\R^{d\times d}_{\rm sym})
\]
such that
\begin{equation}
    \widehat\Sigma_j^\alpha=\mathbb C\ee(\chi_j^\alpha+ P_j^\alpha ) \ \text{in}\ Y\setminus D,\ 
  \diver_y\widehat\Sigma^\alpha_j=0 \ \text{in }\mathcal D'(\T^d;\R^d),
\end{equation}
and
\begin{equation}\label{stress_et_eq1}
    \int_Y \widehat{\Sigma}^\alpha_j dy = \Ahom e^j \otimes e^\alpha, \quad \|\widehat\Sigma^\alpha_j\|_{L^2(Y)}\le C.
\end{equation}
Set
\begin{equation}
  b^\alpha_j:=\widehat\Sigma^\alpha_j-\Ahom e^j \otimes e^\alpha.
\end{equation}
then $b^\alpha_j$ is symmetric, periodic, divergence free, and has zero mean.
\end{proposition}

\begin{proof}
Define the linear map on $H^1(D;\R^d)$:
\[
  \ell^\alpha_j(\phi)
  :=-\int_{Y\setminus D } \mathbb C\ee(\chi^\alpha_j+P_j^\alpha) :\ee ( \mathcal{Q} \phi)\,dy.
\]
For $ r \in \mathcal{R} $, by the definition of $ \mathcal{Q} $ in Lemma \ref{lem:rigid-extension}, $ \mathcal{Q}(r) = \zeta r $, then we have $ e(\mathcal{Q}r) = 0  $ in D, it follows from \eqref{cell_pb} that $ \ell^\alpha_j(r) = 0 $. 

Let
\[
 H^1_\perp(D):=\set{v\in H^1(D ;\R^d):
       \int_D v\cdot r=0\quad\forall r\in\calR}.
\]
Then $ \ell^\alpha_j $ is a bounded linear map on $ H^1_\perp(D) $, and by Korn's inequality, $(u,v)\mapsto \int_D \mathbb C\ee(u):\ee(v)$ is coercive on $H^1_\perp(D)$. By Lax-Milgram, there exists a unique $ v_0\in H^1_\perp(D) $, such that
\begin{equation}
    \int_D \mathbb C\ee(u):\ee(v_0) = \ell^\alpha_j (u), \ \forall u \in H^1_\perp(D),
\end{equation}
and $ \|e(v_0)\|_{L^2(D)} \le C \| \mathbb C\ee(\chi^\alpha_j +P_j^\alpha ) \|_{L^2(Y)} \le C $.

Let
\begin{equation}
    \widehat{\Sigma}^\alpha_j = \left\{
    \begin{array}{ll}
       \mathbb C\ee(\chi^\alpha_j +P_j^\alpha ),  &  \mathrm{in}\ Y\setminus D, \\
        \mathbb C \ee( v_0 ), &  \mathrm{in}\ D,
    \end{array}
    \right.
\end{equation}
then 
\begin{equation}
    \|\widehat\Sigma^\alpha_j\|_{L^2(Y)}\le C\|e(v_0)\|_{L^2(D)} + C \|\ee(\chi^\alpha_j +P_j^\alpha ) \|_{L^2(Y)} \le C.
\end{equation}
For $ \psi \in H^1(\T^d; \R^d) $, denote $ \psi_1 =  \psi - \mathcal{Q}\psi $, then $ \ee(\psi_1) = 0 $ in $ D $, by \eqref{cell_pb}, we then have 
\begin{equation}
    \int_Y \widehat{\Sigma}^\alpha_j :  \ee (\psi_1 ) = \int_{Y\setminus D} \mathbb C \ee(\chi^\alpha_j +P_j^\alpha ) : \ee ( \psi_1   ) =0.
\end{equation}
Hence, 
\begin{equation}
    \begin{aligned}
        \int_Y \widehat{\Sigma}^\alpha_j : \ee (\psi) & = \int_Y \widehat{\Sigma}^\alpha_j : \big( \ee (\psi_1 ) + \ee( \mathcal{Q} \psi ) \big) \\
        &= \int_Y \widehat{\Sigma}^\alpha_j :  \ee( \mathcal{Q} \psi- \mathcal{Q}r_\psi )+ \int_{Y\setminus D} \mathbb C \ee(\chi^\alpha_j +P_j^\alpha ): \ee(\mathcal{Q}r_\psi)  \\
        & = - \ell^\alpha_j(\psi-r_{\psi}) + \int_D \mathbb C \ee(v_0):\ee( \psi - r_\psi ) = 0,
    \end{aligned}
\end{equation}
here $ r_{\psi } $ is the $L^2(D)$-orthogonal projection onto $\calR$, for the second equality we used $ \ee(\mathcal{Q}r_\psi)|_D = \ee(r_\psi)|_D =0 $, and for the third equality we used $ \psi- r_\psi \in H^1_\perp(D) $.

For any $ 1\le i,\beta \le d $, use the method above and the fact that $ \widehat\Sigma^\alpha_i = \widehat\Sigma_\alpha^i $, we have
\begin{equation}
    \begin{aligned}
        \left(\int_Y \widehat\Sigma^\alpha_i \right):e^j\otimes e^\beta
 &=\int_{ Y \setminus D}\mathbb C \ee(\chi^\alpha_i +P_i^\alpha ):\ee(P_j^\beta)  -\int_D\mathbb{C} e(v_0):\ee_y\chi^\beta_j \\
 &=\int_{Y \setminus D }\mathbb C\ee(\chi^\alpha_i +P_i^\alpha ):\ee(\chi^\beta_j +P_j^\beta ) =(\Ahom)_{ij}^{\alpha\beta} ,
    \end{aligned}
\end{equation}
hence the left equality of \eqref{stress_et_eq1} holds.
\end{proof}

We also have the following modification for $ \mathbb{C} \ee(u_\eps) $, which is important in both the proof of qualitative homogenization and the Lipschitz regularity. The proof of it is similar to Proposition \ref{prop:whole-cell-stress}.

\begin{proposition}[Divergence-free filling of the inclusion stresses]\label{prop:solution-stress-fill}
Let $u_\eps\in V_\eps$ solve \eqref{Lame_eq_main}, there exists a symmetric field $\widehat\sigma_\eps\in L^2(\Omega;\R^{d\times d}_{\rm sym}) $, such that
\begin{equation}
    \widehat{\sigma}_\eps = \mathbb{C} \ee(u_\eps) \ \mathrm{in}\ \Omega \setminus D_\eps,  \quad \mathrm{and}\ \norm{\widehat\sigma_\eps}_{L^2(\Omega)}\le C\norm{\ee(u_\eps)}_{L^2(\Omega\setminus D_\eps)},
\end{equation}
and 
\begin{equation}
    \int_\Omega\widehat\sigma_\eps:\ee(\phi) =\pair g\phi_{\partial\Omega}, \ \forall \phi\in H^1(\Omega;\R^d).
\end{equation}
That is,
\[
  \diver\widehat\sigma_\eps=0\quad\text{in }\Omega,
  \qquad
  \widehat\sigma_\eps N=g\quad\text{on }\partial\Omega
\]
in the distributional and normal-trace senses.
\end{proposition}

\begin{proof}
For each $ k\in J_\eps $, denote $ D_\eps^k =\eps(k+D) $, define the linear map:
\begin{equation}
    \ell_{\eps,k}(\phi):=-\int_{\eps( k+ Y\setminus D )}
       \mathbb{C} \ee(u_\eps):\ee( \mathcal{Q}_{\eps,k} (\phi) ),
  \qquad \phi\in H^1(D_\eps^k;\R^d).
\end{equation}
here $ \mathcal{Q}_{\eps,k} \phi =  \mathcal{Q}[\phi (\eps(\cdot + k))](\frac{\cdot}{\eps} - k ) $ with $ \mathcal{Q} $ defined in Lemma \ref{lem:rigid-extension}. By the same argument in Proposition \ref{prop:whole-cell-stress}, $ \ell_{\eps,k}|_{\mathcal{R}} = 0 $, hence $ \ell_{\eps,k} $ is a bounded linear map on 
\begin{equation}
    H^1_\perp(D^k_\eps):= \big\{v\in H^1(D_\eps^k ;\R^d): \int_{D_\eps^k} v\cdot r=0\quad\forall r\in\calR \big\}.
\end{equation}
Moreover, there exists a unique $ v_{\eps,k} \in H^1_\perp(D^k_\eps) $, such that
\begin{equation}
    \int_{D^k_\eps}\mathbb C\ee(v_{\eps,k}):\ee(\phi)
  =\ell_{\eps,k} (\phi), 
  \qquad \forall \phi\in H^1_\perp(D^k_\eps),
\end{equation}
and $ \|\ee(v_{\eps,k})\|_{L^2(D_{\eps,k})} \le C \| \mathbb{C} \ee(u_\eps)\|_{L^2(\eps( k+Y\setminus D))} $. Since both sides vanish on rigid motions, the identity above still holds for every $\phi\in H^1(D^k_\eps;\R^d)$.

Let 
\begin{equation}
    \widehat{\sigma}_\eps = \left\{
    \begin{array}{ll}
       \mathbb{C} \ee(u_\eps),  &  \mathrm{in}\ \Omega\setminus D_\eps, \\
        \mathbb C \ee( v_{\eps,k} ), &  \mathrm{in}\ D_\eps^k,\ k\in J_\eps,
    \end{array}
    \right.
\end{equation}
then $\norm{\widehat\sigma_\eps}_{L^2(\Omega)}\le C\norm{\ee(u_\eps)}_{L^2(\Omega\setminus D_\eps)} $. For each $\phi\in H^1(\Omega;\R^d)$, let $ \phi_\eps := \phi - \sum_{k\in J_\eps} \mathcal{Q}_{\eps,k} \phi $, then $ \phi_\eps = 0 $ in $ D_\eps $, and $ \phi_\eps\in\widetilde V_\eps $, hence,
\begin{equation}
    \begin{aligned}
        \int_\Omega\widehat\sigma_\eps:\ee(\phi)& = \int_{\Omega \setminus D_\eps} \mathbb{C} \ee(u_\eps):\ee(\phi_\eps) + \sum_{k\in J_\eps} \int \mathbb C \ee(u_\eps  ): e(\mathcal{Q}_{\eps,k}\phi)\\
        & = \pair g{\phi_\eps}_{\partial\Omega} + \sum_{k\in J_\eps} \bigg( \ell_{\eps,k}(\phi|_{D_{\eps,k }}) +\int_{\eps(k+Y\setminus D)} \mathbb C \ee(u_\eps  ): e(\mathcal{Q}_{\eps,k}\phi)\bigg)  \\
        &= \pair g{\phi_\eps}_{\partial\Omega} .
    \end{aligned}
\end{equation}
By construction, $ \phi = \phi_\eps $ in the neighborhood of $ \partial\Omega $, the proof is completed.
\end{proof}

\subsection{Qualitative homogenization}


We are now ready to prove the qualitative convergence.

\begin{theorem}
        Let $ u_\eps \in V_\eps $ be the solution of \eqref{Lame_eq_main}, and $ u_0 \in H^1_{\mathcal{R}}(\Omega;\R^d) $ be the solution of \eqref{lame_eq_homo}, then 
\begin{equation}
    u_\eps \rightarrow u_0 \ \mathrm{in}\ L^2(\Omega;\R^d), \quad \mathrm{and}\ u_\eps \rightharpoonup u_0\ \mathrm{in}\ H^1(\Omega).
\end{equation}
\end{theorem}

\begin{proof}
    By Korn's Inequality, 
\begin{equation}
    \|u_\eps\|_{H^1(\Omega;\R^d)} \le C \|\ee(u_\eps)\|_{L^2(\Omega;\R^d)} \le  C \|g\|_{H^{-1/2}(\partial\Omega;\R^d)},
\end{equation}
then there exist a subsequence and $ u_0\in H^1(\Omega;\R^d) $, such that
\begin{equation}
    u_\eps \rightarrow u_0 \ \mathrm{in}\ L^2(\Omega;\R^d), \quad \mathrm{and}\ u_\eps \rightharpoonup u_0\ \mathrm{in}\ H^1(\Omega),
\end{equation}
it remains to show that $ u_0 $ is the solution of \eqref{lame_eq_homo}:

By Proposition \ref{prop:solution-stress-fill}, there exists $ \widehat\sigma_\eps \in L^2(\Omega; \R^{d\times d}_{\sym}) $, such that
\begin{equation}
    \widehat{\sigma}_\eps= \mathbb{C}\ee(u_\eps)\ \mathrm{in}\ \Omega \setminus D_\eps, \quad \diver \widehat{\sigma}_\eps = 0 \ \mathrm{ in}\ \Omega , \quad \widehat\sigma_\eps N=g\quad\text{on }\partial\Omega,
\end{equation}
with the energy estimate
\begin{equation}
    \norm{\widehat\sigma_\eps}_{L^2(\Omega)}\le C\norm{\ee(u_\eps)}_{L^2(\Omega\setminus D_\eps)}.
\end{equation}
It follows that there exists $ \sigma_0 \in L^2(\Omega; \R^{d\times d}_{\sym}) $, after extracting a subsequence, 
\begin{equation}
    \widehat \sigma_\eps \rightharpoonup \sigma_0 \quad \mathrm{in}\,L^2(\Omega; \R^{d\times d}_{\sym}),
\end{equation}
since $ \diver\widehat\sigma_\eps = 0 $, we also have $ \diver \sigma_0 =0 $.

\emph{Step 1.} We next show that $ \sigma_0 = \Ahom \ee(u_0) $: 
For each $ 1\le j,\alpha\le d $, take $ \widehat{\Sigma}_{j}^\alpha $ as in Proposition \ref{prop:whole-cell-stress}, take cutoff function $ \zeta\in C_c^\infty (\Omega) $, since $ \ee( \chi_j^\alpha +P_j^\alpha ) = 0 $ in $ D $, $ \ee(u_\eps) = 0 $ in $ D_\eps $, for $ \eps $ small enough, we have:
\begin{equation}
    \int_{\Omega}\zeta\widehat\sigma_\eps: \ee( \chi_j^\alpha +P_j^\alpha )(\frac\cdot\eps) = \int_{\Omega} \zeta\ee(u_\eps): \widehat\Sigma_j^\alpha(\frac\cdot\eps) .
\end{equation}
Integrating by parts, the left hand side becomes
\begin{equation}\label{599}
    \begin{aligned}
        \int_{\Omega}\zeta\widehat\sigma_\eps: \ee( \chi_j^\alpha +P_j^\alpha ) (\frac\cdot\eps) & =  -\int_{\Omega} \widehat{\sigma}_\eps: \nabla\zeta \otimes ( P_j^\alpha + \eps\chi_j^\alpha(\frac{\cdot}{\eps}) )\\ 
        & \rightarrow -\int_{\Omega}\sigma_0: \nabla\zeta \otimes  x_je^\alpha = \int_{\Omega} \zeta \sigma_0 :(e^\alpha\otimes e^j),
    \end{aligned}
\end{equation}
here we used the symmetry of $ \sigma_0, \widehat \sigma_\eps $ and $ \diver \widehat \sigma_\eps = 0 $, $ \diver \sigma_0 = 0 $ for the first and last equality.

Again by integrating by parts, the right hand side becomes
\begin{equation}
    \begin{aligned}
        \int_{\Omega} \zeta\ee(u_\eps): \widehat\Sigma_j^\alpha(\frac\cdot\eps) & =- \int_{\Omega} \sym\{\nabla\zeta\otimes u_\eps \} : \widehat \Sigma_j^\alpha(\frac\cdot\eps),
    \end{aligned}
\end{equation}
since $ \int_{Y} \widehat\Sigma_j^\alpha = \Ahom (e^j \otimes e^\alpha ) $, we have $ \widehat\Sigma_j^\alpha(\frac{\cdot}{\eps})\rightharpoonup \Ahom (e^j \otimes e^\alpha ) $ in $ L^2(\Omega) $, combine with $ u_\eps \rightarrow u_0 $ in $ L^2(\Omega) $, we have:
\begin{equation}\label{614}
    \int_{\Omega} \zeta\ee(u_\eps): \widehat\Sigma_j^\alpha(\frac\cdot\eps) \rightarrow -\int_\Omega \sym\{\nabla\zeta\otimes u_0 \} : \Ahom (e^\alpha\otimes e^j ) = \int_{\Omega} \zeta \ee(u_0) : \Ahom ( e^\alpha\otimes e^j ).
\end{equation}

Compare the right hand side of \eqref{599} and \eqref{614}, we then derive:
\begin{equation}
    \sigma_0 = \Ahom \ee(u_0).
\end{equation}

\emph{Step 2.} We will show that $ u_0 \in H^1_{\mathcal{R}}(\Omega;\R^d) $, and is the solution of \eqref{lame_eq_homo}.
For any $ \phi \in H^1(\Omega;\R^d) $, by the construction of $ \widehat\sigma_\eps $, we have:
\begin{equation}
    \int_{\Omega} \widehat\sigma_\eps : \ee(\phi) = \pair {g }{\phi}_{\partial\Omega}  ,
\end{equation}
let $ \eps $ tend to $0$, we then derive:
\begin{equation}
    \int_{\Omega} \Ahom \ee(u_0) : \ee(\phi) = \pair {g }{\phi}_{\partial\Omega},\quad \forall \phi\in H^1(\Omega;\R^d),
\end{equation}
that is, $u_0 $ solves \eqref{lame_eq_homo}. Since $ u_\eps \rightharpoonup u_0 $ in $ H^1(\Omega) $, we have also $  u_\eps \rightharpoonup u_0 $ in $ H^{1/2}(\partial \Omega) $, hence
\begin{equation}
    \int_{\partial \Omega} u_0 \cdot r = \lim_{\eps\rightarrow0} \int_{\partial \Omega} u_\eps \cdot r = 0, \quad \forall r\in \calR,
\end{equation}
that is, $ u_0 \in H^1_{\mathcal{R}}(\Omega;\R^d) $.

Since every subsequence has a further subsequence converging to the same unique limit $u_0$, we have in fact the whole sequence converges.
\end{proof}

\subsection{Preliminaries}

We collect some basic notations and lemmas that are frequently used in the quantitative part of the argument. Throughout the paper we assume that $\varepsilon \in (0,1)$.

The boundary cut-off function $\eta_\eps \in C_0^{\infty} (\Omega; [0,1])$ is chosen such that
\begin{equation}\label{cut-off_def}
	\left\{
	\begin{aligned}
		& \mathrm{supp}\,\eta_{\varepsilon}\subset \Omega \setminus \mathcal{O}_{(d+1)\varepsilon}, \\
		&\eta_{\varepsilon} \equiv 1 \quad \mathrm{in}\  \Omega \setminus \mathcal{O}_{(d+2)\varepsilon}, \\
		& \|\nabla \eta_{\varepsilon}\|_{L^{\infty}(\Omega)} \leq C\varepsilon^{-1} ,
	\end{aligned}\right.
\end{equation}
where $\mathcal{O}_{r}$, $r>0$, denotes the boundary layer $\{x\in \Omega: \mathrm{dist}\,(x,\partial \Omega) < r\}$.

Fix a function $\xi \in C_0^{\infty}(B_{1/2})$ such that $0\leq \xi (x)\leq 1$ and $\int_{\mathbb{R}^d}\xi (x)\,dx=1$. We define the smoothing operator $S_{\varepsilon}$ by
\begin{equation}\label{epsi_smoother}
	S_{\varepsilon}(v)(x):=\varepsilon^{-d} \int_{\mathbb{R}^d} v(x-z)\xi(z/\varepsilon)\,dz, \qquad \forall v\in L^1_{\rm loc}(\mathbb{R}^d).
\end{equation}


\begin{lemma}\label{Svarepsilonproperty1}
    Let $U\subset \mathbb{R}^d$ be an open set. For any $ f \in L^2(Y)$, let $\tilde{f}: \mathbb{R}^d \rightarrow \mathbb{C}$ be the periodic extension of $f$. There is a universal constant $C>0$ such that for every $f \in L^2(\mathbb{T}^d)$ and $v \in L^2(\mathbb{R}^d)$,
    \begin{equation}
        \| \tilde{f}^{\varepsilon} S_\varepsilon(v)\|_{L^2( U)}\leq C \|f \|_{ L^2(Y ) } \|v \|_{L^2(U+B_\varepsilon)}.
    \end{equation}
\end{lemma}
\begin{proof}
See \cite[Proposition 3.1.5]{shen2018periodic}. 
\end{proof}

\begin{lemma}
\label{propertySepsilondifference}
    Let $U\subset \mathbb{R}^d$ be an open set. There exists $C=C(d)$ such that for every $v\in H^1(\mathbb{R}^d)$, 
    \begin{equation}
        \left\|S_{\varepsilon}( v )-v  \right\|_{L^2( U)}\leq C \varepsilon  \| \nabla v \|_{L^2(U+B_\varepsilon) } .
    \end{equation}
\end{lemma}
\begin{proof}
By the density argument, we may assume that $ v \in C_0^1(\mathbb{R}^d)$. We compute
\begin{equation}
    \begin{aligned}
        S_{\varepsilon}( v )(x)- v(x)  & = \int_{ |y|\le \frac{1}{2} } \big( v(x-\varepsilon y) -v(x) \big) \xi(y)  \,dy\\
    & =- \int_{ |y|\le \frac{1}{2} } \left(\int_0^1  \varepsilon y \cdot \nabla v(x-t\varepsilon y) \,dt \right) \xi(y)\,dy \\
    & = -\int_{ |y|\le \frac{1}{2} }  \int_0^1   \varepsilon y \cdot\nabla v(x-t\varepsilon y) \xi(y)\,dt  dy.
    \end{aligned}
\end{equation}
It follows that 
\begin{equation}
    \begin{aligned}
        \left\|S_{\varepsilon}( v )-v  \right\|_{L^2( U)}   &\leq  \int_{ |y|\le \frac{1}{2} }  \int_0^1   \| \varepsilon y \cdot\nabla v(\cdot-t\varepsilon y) \xi(y) \|_{L^2(U)}\,dt  dy  \\ &\leq C \varepsilon  \| \nabla v \|_{L^2(U+B_\varepsilon) } .
    \end{aligned}
\end{equation}
\end{proof}

\begin{lemma}\label{blelem}
For any  $u\in H^1(\Omega)$, we have the boundary layer estimate:
\begin{equation}
    \| u \|_{L^2(\mathcal{O}_\varepsilon)} \leq C \varepsilon^{1/2} \|u \|_{H^1(\Omega)},
\end{equation}
where the constant $C>0$ depends only on $\Omega$.
\end{lemma}
\begin{proof}
    See \cite[Proposition 3.1.7]{shen2018periodic}.
\end{proof}



The flux corrector plays an important role in homogenization theory.

\begin{lemma}
	\label{thmflux}
	Suppose that $b=(b_{i}^j )_{i,j=1}^d $ is a tensor field on $Y$ satisfying
	\begin{equation}\label{requireflux}
		\frac{\partial }{\partial y_i} b_{i}^j (y) =0, \qquad \mathrm{for\ all\ } y\in Y,\ 1\leq j \leq d,  
	\end{equation}
	and
	\begin{equation}
		\int_Y b_{i}^j(y)\,dy =0, \qquad \mathrm{for\ all\ } 1\leq i,j \leq d.
	\end{equation}
	Then there exists a tensor field $\Psi = (\Psi_{kij} )_{k,i,j=1}^d$ on $\mathbb{T}^d$ such that
	\begin{equation}
		\label{eq:FtoPsi}
		\frac{\partial }{\partial y_k} \Psi_{kij} (y)= b_{i}^j (y), \qquad \mathrm{for\ all\ } y\in Y, \ 1\leq i,j\leq d,
	\end{equation}
	and
	\begin{equation}	\label{eq:FtoPsi2}
		\Psi_{kij} (y)= - \Psi_{ikj} (y), \qquad \mathrm{for\ all\ } y\in Y,\  1\leq k,i,j\leq d.
	\end{equation}
Moreover, there exists $C=C(d)$ such that 
	\begin{equation}\label{Psib}
		\| \Psi \|_{H^1(\mathbb{T}^d)} \leq C  \|b\|_{L^2(Y)} .
	\end{equation}
	We refer to $\Psi$ as the flux corrector associated with $b$. 
\end{lemma}
\begin{proof}
See \cite[Proposition 3.1.1]{shen2018periodic}.
\end{proof}

\subsection{$ H^1 $-convergence rate}


To apply the usual energy estimate, we have to modify $ u_0 + \eps \chi^\alpha_i (x/\eps) S_\eps (\eta_\eps \partial_i u_0^\alpha) $ into an element of $ V_\eps $. The global convergence rate and the local $ \Ahom $-harmonic approximation used in the Lipschitz regularity below rely on the same first-order calculation. We will first construct the rigidification on $D_\eps$.

The only additional ingredient for a localization domain is a correction in the cells cut by its artificial boundary, the construction of which we will postpone to Section \ref{section_caccioppoli}.

For $u\in H^1(\Omega;\R^d)\cap H^2_{\rm loc}(\Omega;\R^d)$, define its boundary-layer energy as:
\begin{equation}\label{eq:B-eps-def}
  \mathcal{E}(u)
  :=\norm{\ee(u)}_{L^2( \mathcal{O}_{(d+2)\eps})}
    +\eps\norm{\nabla\ee(u)}_{L^2(\Omega\setminus \mathcal{O}_{(d+1)\eps})}.
\end{equation}
here $ \mathcal{O}_\eps := \{x\in \Omega : \mathrm{dist} (x,\partial \Omega) <\eps \}  $ is the boundary layer.


\begin{proposition}
    \label{prop:complete-cell-consistency}
    For each $v\in H^1(\Omega;\R^d)\cap H^2_{\rm loc}(\Omega;\R^d)$, there exists $ v_\eps\in \widetilde  V_\eps $ such that for each $\varphi\in \widetilde V_\eps$,
\begin{equation}
 \left| \int_\Omega \big(\Ahom\ee(v)- \mathbb C \ee(v_\eps) \big):\ee(\varphi) \right| \le C\mathcal{E} (v)\norm\varphi_{H^1(\Omega)},
\end{equation}
moreover, $ v|_{\partial \Omega} = v_\eps|_{\partial \Omega} $ in trace sense, and
    \begin{equation}
  \norm{\nabla( v+ \eps \chi^\eps S_\eps(\eta_\eps \ee(v) )- v_\eps)}_{L^2(\Omega)}
  +\eps^{-1}\norm{ v+ \eps \chi^\eps S_\eps(\eta_\eps \ee(v) )- v_\eps}_{L^2(\Omega)}
  \le C\mathcal{E}(v).
    \end{equation}
    Here we write
    \begin{equation}
       \chi_i^\alpha(x/\eps) S_\eps (\eta_\eps \partial_i  v^\alpha) = \chi^\eps S_\eps (\eta_\eps \ee(v) ).
    \end{equation}
    for short.
\end{proposition}

\begin{proof}
Denote $ w_\eps = v+ \eps \chi^\eps S_\eps(\eta_\eps \ee(v)) $ for short. For $ j\in J_\eps$, let 
\begin{equation}
    v_\eps = \mathcal{P} [w_\eps(\eps(\cdot + j))] (\frac{\cdot}{\eps}-j), \quad \mathrm{for}\ x\in \eps(j+Y),
\end{equation}
and let $ v_\eps =v $ for $ x\in \Omega\setminus \bigcup_{j\in J_\eps } \eps(j+Y) $. It is clear that $ v_\eps \in \widetilde V_\eps $ and $ (v-v_\eps)|_{\partial \Omega } = 0 $. By Lemma \ref{lem:rigidification} and scaling, we have:
\begin{equation}
    \eps^{-1}\|v_\eps -w_\eps\|_{L^2(\bigcup_{j\in J_\eps } \eps(j+Y))} + \|\nabla (v_\eps -w_\eps )\|_{L^2(\bigcup_{j\in J_\eps } \eps(j+Y))} \le C \|e(w_\eps)\|_{L^2(D_\eps )},
\end{equation}
for $ j\notin J_\eps $, since $ S_\eps(\eta_\eps \ee(v))\equiv 0 $ on $ \mathcal{O}_{d\eps} $, $ w_\eps = v = v_\eps $ on $ \eps(j+Y) $. It remains to estimate $\|e(w_\eps)\|_{L^2(D_\eps)}$:
\begin{equation}\label{consruct_Veps_eq1}
    \begin{aligned}
        \| \ee (w_\eps)\|_{L^2(D_\eps)}&  \le \|e(v) - S_\eps(\eta_\eps \ee(v) )\|_{L^2(D_\eps)} + C\eps\|\chi\|_{L^2(Y)}\|\nabla ( S_\eps(\eta_\eps \ee(v) )) \|_{L^2(\Omega)}\\
     &\le \|(1-\eta_\eps) e(v)\|_{L^2(\Omega)} + C \eps \| \nabla(\eta_\eps e(v))\|_{L^2(\Omega)} \le C \mathcal{E} (v).
    \end{aligned}
\end{equation}
here for the first inequality, we used $ e(\chi ^\alpha) + E^\alpha =0  $ in $ D $ and Lemma \ref{Svarepsilonproperty1}, for the second inequality, we used Lemma \ref{propertySepsilondifference}.

We now prove the first statement. For any $\varphi\in \widetilde V_\eps$, by direct computation, we have:
\begin{equation}
    \begin{aligned}
        &\bigg|  \int_{\Omega} (\Ahom e(v) - \mathbb{C} \ee(v_\eps)): \ee(\varphi)\bigg|  \\
        \le& \bigg|\int_{\Omega} \big(\Ahom e^\alpha\otimes e^i - \widehat\Sigma^\alpha_i(x/\eps)\big)S_\eps (\eta_\eps \partial_i v^\alpha): \ee(\varphi) \bigg|  \\
        & + C \big(\|\ee(v) -S_\eps (\eta_\eps \ee(v))\|_{L^2(\Omega)} + \eps  \|\chi\|_{L^2(Y)}\|\nabla ( S_\eps(\eta_\eps \ee(v) )) \|_{L^2(\Omega)} \big)\|\ee(\varphi)\|_{L^2(\Omega)} 
    \end{aligned}
\end{equation}
since $ \mathbb{C} \ee(\chi^\alpha_i+P_i^\alpha) = \widehat\Sigma^\alpha_i (x/\eps) $ in $\mathrm{supp}\, S_\eps (\eta_\eps \ee(v)) \setminus D_\eps $, and $ e(\varphi) = 0 $ in $ D_\eps $. By Proposition \ref{prop:whole-cell-stress}, $ b^\alpha_i = \widehat{\Sigma}_i^\alpha - \Ahom e^\alpha\otimes e^i $ is symmetric, periodic, divergence free, and has zero mean, hence by Lemma \ref{thmflux}, there exists $ \Psi_i^\alpha $ such that $ \partial_l (\Psi_i^\alpha )_{ljk} = (b_i^\alpha)_{j}^k $, then
\begin{equation}
    \begin{aligned}
        & \bigg|\int_{\Omega} \big(\Ahom e^\alpha\otimes e^i - \widehat\Sigma_i^\alpha(x/\eps) \big) S_\eps (\eta_\eps \partial_iv^\alpha ): \ee(\varphi) \bigg| \\
        = & \bigg|\eps \int  \diver (\Psi^\alpha_i(x/\eps)) S_\eps (\eta_\eps \partial_iv^\alpha ): \ee(\varphi) \bigg| \\
         = & \eps \bigg| \int  (\Psi^\alpha_i)_{ljk}(x/\eps)
    \partial_lS_\eps (\eta_\eps \partial_iv^\alpha)\partial_j\varphi^k \bigg|\\
    \le & \eps C \|\Psi^\alpha_i\|_{L^2(Y)} \|\nabla ( S_\eps(\eta_\eps \ee(v) )) \|_{L^2(\Omega)} \|\varphi\|_{H^1(\Omega)},
    \end{aligned}
\end{equation}
here we used the skew-symmetry of $ \Psi^\alpha_i $ and Lemma \ref{Svarepsilonproperty1}.

The first statement then follows from combining the two inequalities above and the same estimate as in \eqref{consruct_Veps_eq1}.
\end{proof}


We are now in the position to prove our main result:

We first give the estimate for the energy $ \mathcal{E}(\cdot) $:
\begin{lemma}[Macroscopic Neumann boundary layer]\label{lem:hom-boundary-layer}
Let $g\in L^2_{\calR}(\partial\Omega)$ and $u_0$ solve \eqref{lame_eq_homo}. Then there exists $ C= C(\Ahom, \Omega,d) $, such that
\begin{equation}\label{eq:hom-boundary-layer}
 \norm{\nabla u_0}_{L^2(\mathcal{O}_t)}
 +t\norm{\nabla^2u_0}_{L^2(\Omega\setminus \mathcal{O}_t)}
 \le Ct^{1/2}\norm g_{L^2(\partial\Omega)}.
\end{equation}
\end{lemma}
\begin{proof}
    See \cite[Lemma 3.3]{geng2023hardy}.
\end{proof}

\begin{proof}[Proof of Theorem~\ref{thm:main-rate}]
Denote $ w_\eps = u_0+ \eps \chi_i^\alpha(x/\eps) S_\eps (\eta_\eps \partial_i u_0^\alpha) $, and $ u_{0,\eps} \in \widetilde V_\eps $ be taken by replacing $ v $ by $ u_0 $ in Proposition  \ref{prop:complete-cell-consistency}. Since $u_{0,\eps} = u_0 $ on the boundary layer of $ \partial\Omega $, we have in fact $ u_{0,\eps} \in V_\eps \subset \widetilde V_\eps \cap H^1_{\calR}(\Omega)$. 

By Proposition  \ref{prop:complete-cell-consistency}, we have: for each $ \varphi\in V_\eps $,
\begin{equation}
    \begin{aligned}
        \bigg|\int_\Omega \mathbb{C} \ee( u_\eps -u_{0,\eps}): \ee(\varphi) \bigg|& = \bigg|\pair g\varphi_{\partial\Omega}-\int_\Omega \mathbb{C} \ee(u_{0,\eps}): \ee(\varphi) \bigg|\\
        & = \bigg|\int_\Omega \big(\Ahom \ee( u_0) - \mathbb{C} \ee(u_{0,\eps})\big): \ee(\varphi) \bigg|\\
         & \le C\mathcal{E} (u_0)\norm\varphi_{H^1(\Omega)}.
    \end{aligned}
\end{equation}
In particular, we can take $ \varphi = u_\eps -u_{0,\eps} $, using Korn's inequality, we have:
\begin{equation}
    \norm {u_\eps - u_{0,\eps}}_{H^1(\Omega)} \le C \mathcal{E} (u_0).
\end{equation}
Again by Proposition  \ref{prop:complete-cell-consistency},
\begin{equation}
    \|u_{0,\eps} - w_\eps\|_{H^1(\Omega)} \le C \mathcal{E} (u_0),
\end{equation}
the Theorem then follows from Lemma \ref{lem:hom-boundary-layer}.
\end{proof}

\section{Caccioppoli-type inequalities and local $\Ahom $-harmonic approximation} \label{section_caccioppoli}

The rest of the article is devoted to proving the global Lipschitz estimate Theorem \ref{thm:main-global}.

We will follow the now standard approximation scheme. The key steps are to establish Caccioppoli inequalities and local $\Ahom $-harmonic solution, here we say $ v\in H^1(B_r) $ is $ \Ahom $-harmonic on $ B_r $ if
\begin{equation}
    \diver \Ahom \ee(v) = 0 \quad \mathrm{in } \ \mathcal{D}'(B_r).
\end{equation}

For the remainder of the paper we assume
\begin{equation}\label{eq:boundary-assumptions}
  \partial\Omega,\partial D\in C^{2,\beta},
  \qquad
  0<\sigma<\beta<1,
  \qquad
  g\in C^{0,\sigma}(\partial\Omega;\R^d)
       \cap H^{-1/2}_{\calR}(\partial\Omega).
\end{equation}

Fix $x_0\in\partial\Omega$.  After a rigid change of coordinates, which we suppress in the notation, there are $r_0>0$ and a function $\psi\in C^{2,\beta}$ such that $x_0=0$, $\psi(0)=0$, $\nabla\psi(0)=0$, 
\begin{align}
  U_r&:=\set{(x',x_d):\abs{x'}<r,
           \ \psi(x')<x_d<\psi(x')+r}\subset\Omega,\\
  \Delta_r&:=\set{(x',\psi(x')):\abs{x'}<r}\subset\partial\Omega.
\end{align}
Let $ \Gamma_r = \partial U_r \setminus \Delta_r $ be the artificial boundary.

In the following, we denote for each ball $ B_r(x_0) \subset\R^d $,
\begin{equation}\label{def_VepsBall}
    V_\eps(B_r(x_0)) = \{u\in H^1(B_r(x_0);\R^d): \ee(u) = 0 \ \mathrm{in}\ B_r(x_0)\cap\eps(D+\Z^d)  \}.
\end{equation}


\subsection{Caccioppoli-type inequalities}


We first construct the adapted cutoff function. The usual cutoff argument must be modified in the presence of exact rigidity. Indeed, multiplying an admissible displacement by a general cutoff destroys the condition \(\ee(u)=0\) inside an inclusion. 

\begin{lemma}\label{lem:inclusion-cutoff} 
For any $0<r<R$, $ R-r>3d\eps $, $x_0\in\R^d$, there exist $\zeta_1\in C_c^\infty(B_R(x_0);[0,1])$ and constant $ C=C(d,D) $ satisfying
\begin{enumerate}[label=\textup{(\roman*)}]
\item $\zeta_1=1$ on $B_r(x_0)$;
\item $\norm{\nabla\zeta_1}_{L^\infty}\le C(R-r)^{-1}$;
\item $\zeta_1$ is constant in a neighborhood of every inclusion $\eps(k+D)$.
\end{enumerate}

Moreover, there exist $\zeta_2\in C^\infty(U_R;[0,1])$ and constant $ C=C(d,D,\Omega) $ satisfying
\begin{enumerate}[label=\textup{(\roman*)},leftmargin=2.5em]
\item  $\zeta_2=1$ on $U_r$, and $\zeta_2=0$ in a neighborhood of $\Gamma_R$;
\item $\abs{\nabla\zeta_2 }\le C(R-r)^{-1}$;
\item $\zeta_2$ is constant in a neighborhood of $\eps(k+D)$, for each $ k\in J_\eps $.
\end{enumerate}
\end{lemma}

\begin{proof}
Denote $r_j = r+ \frac{j}{3}(R-r) $, $ j=1,2 $, for short. Let
\begin{equation}
    \widetilde{D}_j:=\big\{ x\in Y: \dist (x,D) <\frac{j}{4} \delta_0 \big\}, \quad j=1,2,
\end{equation}
take $ \phi \in C_c^\infty(\widetilde D_{2};[0,1]) $,  such that $ \phi \equiv 1 $ in $ \widetilde D_1 $, let $ \phi_\eps^k(x):=\phi(\frac{x}{\eps}-k) $.

For the first statement, take $ \zeta_0 \in C_c^\infty(B_{r_2}(x_0);[0,1]) $ such that $ \zeta_0 \equiv 1 $ on $ B_{r_1} (x_0) $, and $ \| \nabla \zeta_0 \|_{L^\infty (B_{r_2}(x_0)) } \le  C/(R-r) $. Let
\begin{equation}
    \zeta_1 (x):= \big( 1-\sum_{k\in \Z^d} \phi_\eps^k(x)\big) \zeta_0 (x) + \sum_{k\in \Z^d}\phi_\eps^k(x) \zeta_0 (\eps k).
\end{equation}
It is clear that $ 0\le \zeta_1 \le 1 $, and for each $k $, $ \zeta_1 \equiv \zeta_0(\eps  k) $ in $ \eps(k+\widetilde D_1) $. Since $ R-r> 3d\eps $, $ \zeta_1\equiv 1 $ on $ B_{r_1-d\eps}(x_0) \supset B_{r}(x_0) $. By direct computation, for each $ k $, and each $ x\in \eps(k+Y) $,
\begin{equation}
    \begin{aligned}
        \big|\nabla \zeta_1 (x)\big| & =  \big|\nabla\phi_\eps^k(x) (\zeta_0 (x) -\zeta_0 (\eps k)) + (1-\phi_\eps^k(x)) \nabla \zeta_0 (x) \big|\\
        & \le \frac{|x-\eps k|}{\eps} \|\nabla \phi\|_{L^\infty} \|\nabla\zeta_0\|_{L^\infty}  + \|\nabla\zeta_0\|_{L^\infty} \le \frac{C}{R-r}.
    \end{aligned}
\end{equation}

For the second statement,  take $ \zeta \in C^\infty(U_{r_2}(x_0);[0,1]) $ such that $ \zeta \equiv 1 $ on $ U_{r_1} (x_0) $, $ \supp{\zeta} \subset \overline U_{r_2} $, $ \zeta $ vanishes on $ \Gamma_{r_2} $, and $ \| \nabla \zeta \|_{L^\infty (B_{r_2}(x_0)) } \le  C/(R-r) $, here the constant $ C $ only depends on $ \partial\Omega $. Let
\begin{equation}
    \zeta_2 (x):= \big( 1-\sum_{k\in J_\eps } \phi_\eps^k(x)\big) \zeta (x) + \sum_{k\in J_\eps }\phi_\eps^k(x) \zeta(\eps k).
\end{equation}
It is clear that $ 0\le \zeta_2 \le 1 $, for each $k \in J_\eps$, $ \zeta_2 \equiv \zeta(\eps  k) $ in $ \eps(k+\widetilde D_1) $ and $ \zeta_2 \equiv 1 $ on $ U_{r_1-d\eps} \supset U_{r} $. And for each $ k\in J_\eps $, by the same estimate as above, we have:
\begin{equation}
    \sup_{x\in \eps(k+Y)} |\nabla \zeta_2(x) | \le \frac{C}{R-r}.
\end{equation}
For $ x\in U_R \setminus \bigcup_{k\in J_\eps} \eps(k+Y) $, $ \nabla\zeta_2 = \nabla \zeta $, hence the second property holds.
\end{proof}

\begin{theorem}[Caccioppoli inequality]\label{thm:interior-caccioppoli}
Let $u_\eps \in V_\eps (B_R(x_0))$ satisfying
\begin{equation}\label{1108}
    \int_{B_R (x_0)} \mathbb{C} \ee(u_\eps) : \ee(\varphi) = 0, \quad \forall \varphi \in V_\eps(B_R(x_0))\cap  H_0^1(B_R(x_0);\R^d),
\end{equation}
here $ V_\eps (B_R(x_0)) $ is defined by \eqref{def_VepsBall}. Then for $0<r<R$ and $ R-r\ge 3d\eps$,
\begin{equation}
  \int_{B_r(x_0)}\abs{\ee(u_\eps)}^2\,dx
  \le\frac{C}{(R-r)^2}\inf_{c\in\R^d}
        \int_{B_R(x_0)}\abs{u_\eps-c}^2\,dx.
\end{equation}
here the constant depends only on $ d,\mu,\lambda,D $.
\end{theorem}

\begin{proof}
Fix $c\in\R^d$ and take the cutoff from \cref{lem:inclusion-cutoff}.  Since $\zeta_1$ is constant near every inclusion and $\ee(u_\eps)=0$ in each inclusion, $\varphi=\zeta_1^2(u_\eps-c) \in V_{\eps}(B_R) \cap H^1_0(B_R) $. Testing \eqref{1108}, using Cauchy Schwartz inequality gives
\begin{equation}
    \int_{B_R\setminus \eps(D+\Z^d)}\zeta_1^2\abs{\ee(u_\eps)}^2
  \le C\int_{B_R}\abs{u_\eps-c}^2\abs{\nabla\zeta_1}^2.
\end{equation}
\end{proof}

\begin{lemma}[Scale-invariant mixed Korn and trace estimates]\label{lem:mixed-Korn-trace}
There exists a constant $ C= C(\partial\Omega) $, such that for each $v\in H^1(U_R;\R^d)$ such that $v=0$ on $\Gamma_R$ in the sense of trace, we have:
\begin{equation}
  R^{-1}\norm v_{L^2(U_R)}+\norm{\nabla v}_{L^2(U_R)}+R^{-1/2}\norm v_{L^2(\Delta_R)}
  \le C\norm{\ee(v)}_{L^2(U_R)}
\end{equation}
\end{lemma}
\begin{proof}
    The result follows from Korn's Inequality, Trace Theorem and scaling.
\end{proof}

\begin{proposition}[Boundary Caccioppoli inequality]\label{prop:boundary-caccioppoli}
Let $u_\eps\in V_\eps$ be the solution of \eqref{Lame_eq_main}. Assume $ 0<r<R $, $ R-r>3d\eps $, then
\begin{equation}\label{eq:boundary-caccioppoli}
  \int_{U_r}\abs{\ee(u_\eps)}^2
  \le\frac C{(R-r)^2} \inf_{c\in\R^d}\int_{U_R}\abs{u_\eps-c}^2
     +C R^d \norm g_{L^\infty(\Delta_R)}^2.
\end{equation}
In particular, for $ R> 12d\eps $,
\begin{equation}\label{eq:boundary-caccioppoli-fixed}
  \left(\avgint_{U_{7R/4}}\abs{\ee(u_\eps)}^2\right)^{1/2}
  \le C\left\{
      \frac1R\inf_{c\in\R^d}
         \left(\avgint_{U_{2R}}\abs{u_\eps-c}^2\right)^{1/2}
      +\norm g_{L^\infty(\Delta_{2R})}\right\}.
\end{equation}
\end{proposition}

\begin{proof}
Let $\zeta_2$ be constructed as in \cref{lem:inclusion-cutoff}. Since $\zeta_2$ is constant in the neighborhood of every inclusion and $\ee(u_\eps)=0$ there, $\varphi:=\zeta_2^2(u_\eps-c) \in V_\eps ( U_R )$, and vanishes on $\Gamma_R$. Testing gives
\begin{align}
 \int_{U_R\setminus D_\eps}\zeta_2^2\mathbb C\ee(u_\eps):\ee(u_\eps)
 &=-2\int_{U_R\setminus D_\eps}
      \zeta_2\mathbb C\ee(u_\eps):
      \sym((u_\eps-c)\otimes\nabla\zeta_2)\notag\\
 &\quad+\int_{\Delta_R}g\cdot\zeta_2^2(u_\eps-c).
 \label{eq:boundary-caccioppoli-test}
\end{align}
\cref{lem:mixed-Korn-trace} gives
\begin{equation}
    \begin{aligned}
        \norm {\zeta_2 (u_\eps -c)}_{L^2(\Delta_R)}
  &\le CR^{1/2}\norm{\ee(\zeta_2 (u_\eps -c))}_{L^2(U_R)} \\
  &\le C R^\frac{1}{2} \big( \frac{1}{R-r} \norm{u_\eps -c}_{L^2(U_R)}+ \norm{\zeta_2 \ee(u_\eps)}_{L^2(U_R)}\big),
    \end{aligned}
\end{equation}
By Cauchy Schwartz Inequality and $\abs{\Delta_R}\le CR^{d-1}$, we then derive:
\begin{equation}
    \begin{aligned}
        \bigg|\int_{U_R\setminus D_\eps}\zeta_2^2\mathbb C\ee(u_\eps):\ee(u_\eps) \bigg| \le & C (R-r)^{-2}\int_{U_R}\abs{u_\eps-c}^2 + R\|g\|^2_{L^2(\Delta_R)}\\
        \le & C (R-r)^{-2}\int_{U_R}\abs{u_\eps-c}^2 + R^d\|g\|^2_{L^\infty(\Delta_R)}.
    \end{aligned}
\end{equation}
\end{proof}



\subsection{Local $ \Ahom $-harmonic approximation}


To obtain the one-step improvement, we need to construct a local $ \Ahom $-harmonic approximation for a local solution. We will prove the following two results. These comparison estimates provide the quantitative input for the excess-decay arguments in the next section.

\begin{proposition}[Boundary homogenized approximation]\label{prop:boundary-comparison}
There exist constants $c_0, C >0$, such that for any $ 0<\eps\le c_{0}R $, let $ u_\eps \in V_\eps $ be the solution of \eqref{Lame_eq_main}, there exists $u_0\in H^1(U_R;\R^d)$ satisfying
\begin{equation}\label{eq:boundary-hom-problem}
  -\diver(\Ahom\ee(u_0))=0\quad\text{in }U_R,
  \qquad
  \Ahom\ee(u_0)N=g\quad\text{on }\Delta_R,
\end{equation}
and
\begin{equation}\label{eq:boundary-comparison}
  \begin{aligned}
      &\left(\avgint_{U_R}\abs{u_\eps-u_0}^2\right)^{1/2}\\
  & \qquad \qquad\le C\left(\frac\eps R\right)^{1/2} R  \bigg\{ \frac{1}{R}\inf_{c\in\R^d}\left(\avgint_{U_{2R}}\abs{u_\eps-c}^2\right)^{1/2}
       +\norm g_{L^\infty(\Delta_{2R})}\bigg\} ,
  \end{aligned}
\end{equation}
here $ c_0,C $ depend only on $ d,D,\lambda,\mu,\Omega $.
\end{proposition}

\begin{proposition}
    \label{thm:local-compactness}
There exists $ \eps_0= \eps_0(d)>0 $, such that for any $ \eps<\eps_0 $, $ x_0 \in \R^d $, if $ u_\eps\in V_\eps(B_1(x_0)) $ satisfying
\begin{equation}
    \int_{B_1 (x_0)} \mathbb{C} \ee(u_\eps) : \ee(\varphi) = 0, \quad \forall \varphi \in V_\eps(B_1(x_0))\cap H_0^1(B_1(x_0);\R^d) ,
\end{equation}
and
\begin{equation}
  \avgint_{B_1(x_0)}\abs{u_\eps}^2\le1,
\end{equation}
here $ V_\eps (B_1(x_0)) $ is defined as \eqref{def_VepsBall}, then there exists an $\Ahom$-harmonic map $v$ in $B_{3/4} $ with
\begin{equation}\label{eq:compactness-approx}
  \norm{u_\eps-v}_{L^2(B_{3/4}(x_0))}\le C \eps^{1/2},
\end{equation}
here the constant $C $ depends only on $ d,D,\lambda,\mu $.
\end{proposition}

Let $ G $ be a connected $ C^{1,1} $ domain, let $ J_\eps $ be an index set such that
\begin{equation}\label{eq:boundary-cell-selection}
  J_\eps^ 0=\set{k\in \Z^d: \eps(k+ \overline{Y})\subset G}
  \subset J_\eps
  \subset\set{k \in \Z^d:\eps(k+ \overline{Y})\cap G\ne\varnothing},
\end{equation}
define:
\begin{equation}
    D^J_\eps : = G \cap \bigcup_{k\in J_\eps} \eps(k+D), \quad \mathrm{and}\ V_{\eps}(G):=\set{v\in H^1(G;\R^d):
       \ee(v)=0\text{  in }D^J_{\eps}}.
\end{equation}

\begin{lemma}
\label{thm:dual-consistency}
For any $ M>1 $, there exists a constant $ C = C (\lambda, \mu,d,D,M) $, such that for any simply connected domain $ G $ with its boundary $ C^{1,1} $-norm less than $ M $, the following holds:
    
    For each $ v\in  H^2(G;\R^d)$, there exists $ v_\eps \in V_\eps (G) $ such that
    \begin{equation}\label{4357}
        \|v-v_\eps\|_{L^2(G)} + \|v-v_\eps\|_{H^{1/2}(\partial G)}
  \le C\eps^{1/2}\norm v_{H^2(G)},
    \end{equation}
    and for each $ \varphi \in V_\eps(G) $, 
    \begin{equation}
 \left|
  \int_{G}\mathbb C\ee(v_\eps):\ee(\varphi)
  -\int_G\Ahom\ee(v):\ee(\varphi)
 \right|
 \le C\eps^{1/2}\norm v_{H^2(G)}\norm\varphi_{H^1(G)}.
\end{equation}
\end{lemma}

\begin{proof}
Since $ v\in  H^2(G;\R^d) $, by Lemma \ref{blelem}, we have:
\begin{equation}
    \norm{\ee(v)}_{L^2( G_{(d+2)\eps})}
    +\eps\norm{\nabla\ee(v)}_{L^2(G\setminus G_{(d+1)\eps})} \le C\varepsilon^{\frac{1}{2}} \|v\|_{H^2(G;\R^d)},
\end{equation}
here $ G_t:=\{x\in G: \dist(x,\partial G)<t\} $ is the boundary layer of $ G $.

Replace $ \Omega $ by $ G $ in Proposition \ref{prop:complete-cell-consistency}, we can construct $ v_{\eps,0} \in H^1 ( G; \R^d ) $ such that $ e(v_{\eps,0}) = 0 $ in $ \bigcup_{j\in J_\eps^0} \eps(j+D) $, such that
\begin{equation}\label{4380}
    \norm{ v+ \eps \chi^\eps S_\eps(\zeta_\eps \ee(v) )- v_{\eps,0} }_{H^1(G)} \le C\varepsilon^{1/2} \|v\|_{H^2(G;\R^d)},
\end{equation}
and for any $ \psi \in H^1(G;\R^d) $, such that $ e(\psi) = 0 $ in $ \bigcup_{j\in J_\eps^0} \eps(j+D) $,
\begin{equation}\label{4384}
    \left| \int_G \big(\Ahom\ee(v)- \mathbb C \ee(v_{\eps,0}) \big):\ee(\psi) \right| \le C\varepsilon^{1/2} \|v\|_{H^2(G;\R^d)}\norm\psi_{H^1(G )},
\end{equation}
here $ \zeta_\eps\in C_c^\infty (G;[0,1]) $ is a cutoff function such that $ \zeta_\eps \equiv 1 $ in $ G\setminus G_{(d+2)\eps} $, $ \supp \zeta_\eps \subset G\setminus G_{(d+1)\eps} $, $ \|\nabla \zeta_\eps \|_{L^\infty ( G )} \le C\eps^{-1} $.

Take an extension operator $ E_G :H^2(G;\R^d) \rightarrow H^2(\R^d;\R^d)$, such that $ \|E_G \|_{H^2(G) \rightarrow H^2(\R^d)} \le C $. Let $ v_\eps = v_{\eps,0} $ in $ G \setminus \bigcup_{j\in J_\eps\setminus J_\eps^0}\eps(j+Y) $, and for each $ j\in J_\eps\setminus J_\eps^0 $, take
\begin{equation}
    v_\eps = \mathcal{P} \big[(E_G v)(\eps(\cdot + j))\big] (\frac{\cdot}{\eps}-j), \quad \mathrm{for}\ x\in \eps(j+Y)\cap G,
\end{equation}

It is easy to check that $ v_\eps \in V_\eps(G) $, and by Lemma \ref{lem:rigidification} and scaling, for $ j\in J_\eps \setminus J_\eps^0 $, $ \|v_\eps - v\|_{H^1(\eps(j+Y)\cap G)}\le C \|e(E_G v)\|_{L^2(\eps(j+D))} $, hence, by \eqref{4380},
\begin{equation}
    \begin{aligned}
        \|v_\eps - (v+ \eps \chi^\eps S_\eps(\zeta_\eps \ee(v) ))\|_{H^1(G)} 
        \le & \|v_{\eps,0} - (v+ \eps \chi^\eps S_\eps(\zeta_\eps \ee(v) ))\|_{H^1(\bigcup_{j\in J_\eps^0} \eps(j+Y))} \\
        &+ \| v_\eps - v \|_{H^1(\bigcup_{j\in J_\eps\setminus J_\eps^0} \eps(j+Y))}\\
        \le &  C\varepsilon^{1/2} \|v\|_{H^2(G;\R^d)} + C \|e(E_G v)\|_{L^2(G_{d\eps})} \\
        \le & C\varepsilon^{1/2} \|v\|_{H^2(G;\R^d)},
    \end{aligned}
\end{equation}
here we used $ \bigcup_{j\in J_\eps\setminus J_\eps^0} \eps(j+Y) \subset G_{d\eps} $, and $ \supp S_\eps(\zeta_\eps \ee(v) ) \subset G\setminus G_{d\eps} $. \eqref{4357} then follows from the Trace Theorem
\begin{equation}
    \|v_\eps - v\|_{H^{1/2} (\partial G) } \le C \|v_\eps - (v+ \eps \chi^\eps S_\eps(\zeta_\eps \ee(v) ))\|_{H^1(G)} \le C\varepsilon^{1/2} \|v\|_{H^2(G;\R^d)},
\end{equation}
and $ \eps\|\chi^\eps S_\eps(\zeta_\eps \ee(v) )\|_{L^2(G)} \le C \eps \|\chi\|_{L^2(Y)} \|v\|_{H^1(G;\R^d}) $.

For the second statement, since for each $ \varphi \in V_\eps(G) $, we have $ e(\varphi) = 0 $ in $ \bigcup_{j\in J_\eps^0} \eps(j+D) $, then \eqref{4384} gives: 
\begin{equation}
    \begin{aligned}
        &\left| \int_G \big(\Ahom\ee(v)- \mathbb C \ee(v_{\eps}) \big):\ee(\varphi) \right|\\ \le & \left| \int_{G}\mathbb C\ee(v_\eps - v_{\eps,0}):\ee(\varphi)\right|+\left| \int_G \big(\Ahom\ee(v)- \mathbb C \ee(v_{\eps,0}) \big):\ee(\varphi) \right| \\
 \le & C \|v_\eps - v_{\eps,0}\|_{H^1(G)} \| e(\varphi)\|_{L^2(G)} +C\eps^{1/2}\norm v_{H^2(G)}\norm\varphi_{H^1(G)}\\
 \le & C\eps^{1/2}\norm v_{H^2(G)}\norm\varphi_{H^1(G)},
    \end{aligned}
\end{equation}
here for the last inequality, we used
\begin{equation}
     \begin{aligned}
         \|v_\eps - v_{\eps,0}\|_{H^1(G)} &= \|v_\eps - v\|_{H^1(G\setminus \bigcup_{j\in J_\eps^0}\eps(j+Y))} \\
         &\le \|v_\eps - (v+ \eps \chi^\eps S_\eps(\zeta_\eps \ee(v) ))\|_{H^1(G)} \le C\varepsilon^{1/2} \|v\|_{H^2(G;\R^d)}.
     \end{aligned}
\end{equation}
\end{proof}

\begin{proposition}
    \label{thm_neumann_approx}
    For any $ M>1 $, there exists a constant $ C = C (\lambda, \mu,d,D,M) $, such that for any simply connected domain $ G $ with its boundary $ C^{1,1} $-norm less than $ M $, the following holds:
    
    Let $u_\eps\in V_\eps(G)$ be a solution of
\begin{equation}
  \int_{G}\mathbb C\ee(u_\eps):\ee(\varphi)
  =\pair h\varphi_{\partial G},
  \qquad  \forall\varphi\in V_\eps (G),
\end{equation}
and let $u_{0}\in H^1(G;\R^d)$ be a solution of
\begin{equation}
  \int_G\Ahom\ee(u_{0}):\ee(\varphi)
  =\pair h\varphi_{\partial G},
  \qquad \forall \varphi\in H^1(G;\R^d).
\end{equation}
here $ h\in H^{-1/2}_{\mathcal{R}} (\partial G;\R^d) $. Then, there exists a rigid motion $r_\eps\in\calR$ such that
\begin{equation}
  \norm{u_\eps-r_\eps-u_{0}}_{L^2(G)}
  \le C\eps^{1/2}\norm h_{H^{-1/2}(\partial G)}.
\end{equation}
\end{proposition}
\begin{proof}
Without loss of generality, assume the $L^2(G)$-orthogonal projection of $ u_0 $ onto $\calR$ is $ 0 $. Let $ r_\eps $ be the $L^2(G)$-orthogonal projection of $ u_\eps $ onto $\calR$, then by Korn's inequality,
\begin{equation}
  \norm{u_\eps -r_\eps }_{H^1(G)}+\norm{u_{0}}_{H^1(G)}
  \le C\norm h_{H^{-1/2}(\partial G)}.
\end{equation}
Take $ z\in H^1 (G;\R^d) $, with its $L^2(G)$-orthogonal projection onto $\calR$ being $ 0 $, solves
\begin{equation}
    \int_G\Ahom\ee(z):\ee(\phi)=\int_G ( u_\eps -u_0 -r_\eps)\cdot\phi,
  \qquad\forall\phi\in H^1(G;\R^d),
\end{equation}
Since $ \partial G \in C^{1,1} $, by standard elliptic estimate, we have:
\begin{equation}
    \|z\|_{H^2(G) } \le C \|u_\eps -u_0 -r_\eps\|_{L^2(G)} .
\end{equation}

Let $ z_\eps \in V_\eps (G) $ be constructed from $ z $ as in Lemma \ref{thm:dual-consistency}, then
\begin{equation}
    \begin{aligned}
        &\|u_\eps -u_0 -r_\eps\|_{L^2(G)}^2 \\
        = &\int_G\Ahom\ee(z):\ee(u_\eps -u_0 -r_\eps)\\
        =& \int_G\Ahom\ee(z):\ee(u_\eps -r_\eps) - \pair h z_{\partial G} + \pair h {z_\eps}_{\partial G} - \int_{G}\mathbb C\ee(u_\eps-r_\eps):\ee(z_\eps)\\
          \le & \bigg|\int_G \big(\Ahom\ee(z)- \mathbb C \ee(z_{\eps}) \big):\ee(u_\eps -r_\eps) \bigg| + \big| \pair h {z_\eps -z}_{\partial G}  \big|\\
          \le & C\eps^{1/2}\norm z_{H^2(G)}\norm{u_\eps-r_\eps}_{H^1(G)} + \|z-z_\eps\|_{H^{\frac{1}{2}}(\partial G)}\|h\|_{H^{-1/2}(\partial G)}\\
          \le & C\eps^{1/2}\|u_\eps -u_0 -r_\eps\|_{L^2(G)} \|h\|_{H^{-1/2}(\partial G)}.
    \end{aligned}
\end{equation}
\end{proof}

We are now ready to prove Proposition \ref{prop:boundary-comparison} and Proposition \ref{thm:local-compactness}.

\begin{proof}[Proof of Proposition~\ref{prop:boundary-comparison}]
Take a connected $ C^{1,1} $ domain $ U_R \subset G_R \subset U_{3R/2} $, such that the scaled domain $ \frac{1}{R}\partial G_R $ has a $ C^{1,1} $-norm  less than $ M $, for some $ M>1 $ independent of $ R $.

Let $ h_R:=\widehat\sigma_\eps N
  \in H^{-1/2}(\partial G_R;\R^d)$, here $ N $ is the outward unit normal vector to $ \partial G_R $. Then $ h_R= g  $ on $ \partial G_R \cap\partial \Omega $, and by Proposition \ref{prop:solution-stress-fill}, for each $ r\in\mathcal{R} $, 
  \begin{equation}
      \pair{h_R}r_{\partial G_R}
  =\int_{G_R}\widehat\sigma_\eps:\ee(r)=0,
  \end{equation}
 moreover, by the construction of $ \widehat{\sigma}_\eps $, for each $ \varphi \in V_\eps (G_R) $, 
\begin{equation}
    \int_{G_R}\mathbb C\ee(u_\eps):\ee(\varphi)
  =\pair {h_R} \varphi_{\partial G_R},
  \qquad  \forall\varphi\in V_\eps (G_R),
\end{equation}
Let $ u_0 $ be the solution of 
\begin{equation}
  \int_{G_R} \Ahom\ee(u_{0}):\ee(\varphi)
  =\pair {h_R}\varphi_{\partial G_R},
  \qquad \forall \varphi\in H^1(G_R;\R^d),
\end{equation}
then by Proposition \ref{thm_neumann_approx} and scaling, we have:
\begin{equation}\label{eq:robust-comparison}
  \norm{u_\eps-r_\eps-u_{0}}_{L^2(G_R)} 
  \le C(\eps R )^{1/2}\norm {h_R}_{H^{-1/2}(\partial G_R)}\le C (\eps R )^{1/2}\|\widehat{\sigma}_\eps \|_{L^2(G_R)}.
\end{equation}

Take $ c_{0}< \frac{1}{4d} $, \eqref{eq:boundary-comparison} then follows from the construction of $ \widehat{\sigma}_\eps $ and Proposition \ref{prop:boundary-caccioppoli}
\begin{equation}
    \begin{aligned}
        \|\widehat{\sigma}_\eps \|_{L^2(G_R)} & \le C \|\ee(u_\eps)\|_{L^2(U_{3R/2+ d\eps})} \le C \|\ee(u_\eps)\|_{L^2(U_{7R/4})}\\
        & \le C R^{ d/2- 1}  \left\{
    \inf_{c\in\R^d}
       \left(\avgint_{U_{2R}}\abs{u_\eps-c}^2\right)^{1/2}
    +R\norm g_{L^\infty(\Delta_{2R})}
  \right\}.
    \end{aligned}
\end{equation}
\end{proof}

Proposition \ref{thm:local-compactness} can be proved in a similar way by constructing suitable boundary condition.

\begin{proof}[Proof of Proposition~\ref{thm:local-compactness}]
     Take $ \eps_0 = \frac{1}{8d} $, then $ B_{3/4+d\eps} \subset B_{7/8} $. Let $ \widetilde \sigma_\eps $ be constructed similarly as in Proposition \ref{prop:solution-stress-fill}, such that
    \begin{equation}
        \begin{aligned}
            &\widetilde{\sigma}_\eps = \mathbb{C} \ee(u_\eps) \ \mathrm{in}\ B_{3/4}(x_0) \setminus \eps (D+\Z^d),  \\
            & \norm{\widetilde\sigma_\eps}_{L^2(B_{3/4}(x_0))}\le C\norm{\ee(u_\eps)}_{L^2(B_{7/8}(x_0) \setminus \eps( D+\Z^d))}.
        \end{aligned}
    \end{equation}
Then we have: $ u_\eps \in V_\eps (B_1(x_0)) $ solves:
    \begin{equation}
        \int_{B_{3/4}(x_0)}\mathbb C\ee(u_\eps):\ee(\varphi)
  =\pair {\widetilde\sigma_\eps N}\varphi_{\partial B_{3/4}(x_0)},
  \qquad  \forall\varphi\in V_\eps (B_{3/4}(x_0)).
    \end{equation}
   
Let $ u_0 $ be a solution of
\begin{equation}
    \int_{B_{3/4}(x_0)} \Ahom \ee( u_0 ):\ee(\varphi) =\pair {\widetilde{\sigma}_\eps N}\varphi_{\partial B_{3/4}(x_0)}, \ \forall \varphi\in H^1(B_{3/4}(x_0);\R^d),
\end{equation}
Since for each $ r\in\mathcal{R} $, $ \pair {\widetilde \sigma_\eps N}r_{\partial B_{3/4}} = 0 $, it then follows from Proposition \ref{thm_neumann_approx} that there exists a rigid motion $r_\eps\in\calR$ such that
\begin{equation}
  \norm{u_\eps-r_\eps-u_{0}}_{L^2(B_{3/4}(x_0))}
  \le C\eps^{1/2}\norm {\widetilde \sigma_\eps N}_{H^{-1/2}(\partial B_{3/4}(x_0))}.
\end{equation}
\eqref{eq:compactness-approx} then follows from
\begin{equation}
    \begin{aligned}
        \norm {\widetilde \sigma_\eps N}_{H^{-1/2}(\partial B_{3/4}(x_0))} & \le C \norm{\ee(u_\eps)}_{L^2(B_{7/8}(x_0))} \\
        &\le C \bigg( \inf_{c\in\R^d}
        \int_{B_1(x_0)}\abs{u_\eps-c}^2\,dx \bigg)^{1/2} \le C,
    \end{aligned}
\end{equation}
and take $ v= u_0+r_\eps $, here we used the Caccioppoli inequality Theorem \ref{thm:interior-caccioppoli}.
\end{proof}

\section{Lipschitz estimate}

The main task of this section is to prove the following two results:

\begin{theorem}
\label{thm_interior_lip_est}
   Assume $ \partial D \in C^{2,\beta} $. Let $u_\eps\in V_\eps$ be the solution of \eqref{Lame_eq_main}. If
\[
  B_{2R}(x_0)\Subset\Omega,
  \qquad
  \dist(B_{2R}(x_0),\partial\Omega)>d\eps,
\]
then
\begin{equation}\label{eq:global-interior-ball}
  \norm{\nabla u_\eps}_{L^\infty(B_R(x_0))}
  \le\frac CR\inf_{c\in\R^d}
       \left(\avgint_{B_{2R}(x_0)}\abs{u_\eps-c}^2\,dx\right)^{1/2}.
\end{equation}
\end{theorem}

\begin{theorem}[Uniform local boundary Lipschitz estimate]\label{thm:local-boundary-Lipschitz}
There are $R_0>0$ and $C<\infty$ such that, for every $x_0\in\partial\Omega$, every $0<R\le R_0$, and every $0<\eps\le1$,
\begin{equation}\label{eq:local-boundary-Lipschitz}
\begin{split}
  \norm{\nabla u_\eps}_{L^\infty(U_{R/2}(x_0))}
  \le C\Bigg\{&\frac1R\inf_{c\in\R^d}
       \left(\avgint_{U_{2R}(x_0)}\abs{u_\eps-c}^2\right)^{1/2}\\
  & \qquad +\norm g_{L^\infty(\Delta_{2R}(x_0))}
   +R^\sigma[g]_{C^{0,\sigma}(\Delta_{2R}(x_0))}\Bigg\}.
\end{split}
\end{equation}
\end{theorem}

Theorem \ref{thm:main-global} can be easily proved by a covering argument provided Theorem \ref{thm_interior_lip_est} and Theorem \ref{thm:local-boundary-Lipschitz} above hold.

The proof of the subsequent elementary lemma is straightforward, we just omit it.

\begin{lemma}[Affine nondegeneracy on microscopic patches]\label{lem:affine-nondegeneracy}
Assume $ \partial D \in C^{2,\beta} $, $ B_{\rho} (x) \subset Y $ and $x\in \partial D $, then there exists a constant $ C= C(d,D) $, such that for every rigid affine map $r(x)=a+Bx$, with $ B $ skew-symmetric, and for every $c\in\R^d$,
\begin{equation}\label{eq:affine-nondegeneracy}
  \rho\abs B+
  \left(\avgint_{B_{\rho}(x)}\abs{r-c}^2\right)^{1/2}
  \le C \left(\avgint_{B_{\rho}(x)\cap D}\abs{r-c}^2\right)^{1/2}.
\end{equation}
\end{lemma}


\begin{proposition}
    \label{thm:micro-interior}
     For fixed $\kappa\ge1$, there exists a constant $ C= C(\kappa,d,\mu,\lambda,D) $, such that for each $y\in\R^d$ and $0<\rho\le\kappa\eps$, assume $ u_\eps \in V_\eps(B_\rho(y)) $ satisfying
\begin{equation}
    \int_{B_\rho(y)} \mathbb{C} \ee(u_\eps) : \ee(\varphi) = 0, \quad \forall \varphi \in V_\eps(B_\rho(y))\cap H_0^1(B_\rho (y);\R^d),
\end{equation}
here $ V_\eps(B_\rho(y)) $ is defined as \eqref{def_VepsBall}, then
\begin{equation}\label{eq:micro-interior-estimate}
  \norm{\nabla u_\eps}_{L^\infty(B_{\rho/2}(y))}
  \le\frac{C_\kappa}{\rho}\inf_{c\in\R^d}
       \left(\avgint_{B_\rho(y)}\abs{u_\eps-c}^2\,dx\right)^{1/2}.
\end{equation}
\end{proposition}

\begin{proof}
Since \eqref{eq:micro-interior-estimate} is scale-invariant, we will assume $ \eps =1 $ without loss of generality.
That is, we will prove that: for any $ u $ solution of
\begin{equation}
    \mathrm{div} (\mathbb{C} \ee (u )) = 0\ \mathrm{in} \ Y\setminus D+\Z^d,\ \ee(u) = 0 \ \mathrm{in} \  D+\Z^d,
\end{equation}
we have: for any $ 2\delta < \kappa $ and $ y\in \R^d $,
\begin{equation}\label{interior_micro_Linfty_eq1}
  \norm{\nabla u}_{L^\infty(B_{\delta}(y))}
  \le\frac{C_\kappa}{\delta}\inf_{c\in\R^d}\left(\avgint_{B_{2\delta}(y)}\abs{u-c}^2\,dx\right)^{1/2}.
\end{equation}

For $ B_{2\delta} (y) \subset Y\setminus D+\Z^d $, \eqref{interior_micro_Linfty_eq1} is just the result of standard interior $ L^\infty  $ estimate of elliptic equation. 

For $ y\in\partial D+ \Z^d $, and $ B_{2\delta}(y) $ belongs to only one cell, then there is $ j\in \Z^d $, such that $ B_{2\delta}(y) \cap ( D+\Z^d) = B_{2\delta}(y) \cap ( j+D)   $ . Denote
\begin{equation}
    u=a+Bx, \ \mathrm{in}\ B_{2\delta}(y) \cap ( D+\Z^d),\ B \,\text{skew-symmetric},
\end{equation}
Then $ v= u-(a+Bx) $ solves
\[
  -\diver(\mathbb C\ee(v))=0\ \text{in } B_{2\delta}(y)\setminus (D+\Z^d),
  \quad
  v=0 \ \text{on }B_{2\delta}(y)\cap\partial D +\Z^d,
\]
it then follows from the standard $ L^\infty  $ estimate of elliptic equation with Dirichlet boundary condition that:
\begin{equation}
     \norm{\nabla v}_{L^\infty(B_{\delta}(y))}
  \le\frac{C_\kappa}{\delta}\inf_{c\in\R^d}\left(\avgint_{B_{2\delta}(y)}\abs{v-c}^2\,dx\right)^{1/2},
\end{equation}
combine with Lemma \ref{lem:affine-nondegeneracy}, we then derive:
\begin{equation}
    \begin{aligned}
        & \norm{\nabla u}_{L^\infty(B_{\delta}(y))} \\
        \le & \frac{C}{\delta}\inf_{c\in\R^d}\left(\avgint_{B_{2\delta}(y)}\abs{v-c}^2\,dx\right)^{1/2}+ |B |\\
        \le & \frac{C}{\delta}\inf_{c\in\R^d}\left(\avgint_{B_{2\delta}(y)}\abs{u-c}^2\,dx\right)^{1/2}+ \frac{C}{\delta} \inf_{c\in\R^d} \left(\avgint_{B_{2\delta}(y)} \abs{r-c}^2\,dx\right)^{1/2} + |B |\\
        \le & \frac{C}{\delta}\inf_{c\in\R^d}\left(\avgint_{B_{2\delta}(y)}\abs{u-c}^2\,dx\right)^{1/2},
    \end{aligned}
\end{equation}
since $ r= u $ in $ B_\delta(y) \cap (D+\Z^d) $.

For $ B_{2\delta} (y) \subset D+\Z^d $, \eqref{interior_micro_Linfty_eq1} is just the result of affine norm equivalence.

For a general $ B_\delta(y) $, \eqref{interior_micro_Linfty_eq1} follows from a covering argument.
\end{proof}

\begin{proposition}
    \label{thm:micro-boundary}
Assume $ \partial \Omega ,\partial D$ are $ C^{2,\beta} $. For fixed $\kappa\ge1$, there is $C= C(\kappa,d,\mu,\lambda,D,\Omega)<\infty$ such that, for each $x_0\in\partial\Omega$, $ 0< \rho \le \kappa \eps  $ and $ u_\eps \in V_\eps $ solution of \eqref{Lame_eq_main},
\begin{equation}\label{eq:micro-boundary-estimate}
\begin{split}
  \norm{\nabla u_\eps}_{L^\infty(U_{\rho/2}(x_0))}
   \le C_\kappa\bigg\{\frac1\rho\inf_{c\in\R^d}
       &\left(\avgint_{U_\rho(x_0)}\abs{u_\eps-c}^2\right)^{1/2} \\
       &\quad +\norm g_{L^\infty(\Delta_\rho(x_0))}
   +\rho^\sigma[g]_{C^{0,\sigma}(\Delta_\rho(x_0))}\bigg\}.
\end{split}
\end{equation}
\end{proposition}

\begin{proof}
Since \eqref{eq:micro-boundary-estimate} is scale-invariant, we will assume $ \rho =1 $ without loss of generality, the condition then turns to be $ \eps >\kappa^{-1} $. 

For fixed $ \eps $, we can always find a finite covering $ \{ B_{\delta_j} (y_j) \}_{j=1}^{N_\eps} $ of $ U_\frac{1}{2}(x_0) $, with $ \delta_j $ comparable to $ \min\{1,\eps\} $, such that $ N_\eps \le N_\kappa $, and for each $ j\in  N_\eps $, we have $ B_{2\delta_j} (y_j) \cap\Omega \subset U_1(x_0) $ and either $ y_j \in \Delta_1 \cup ( U_1(x_0) \cap \partial \eps ( D+\Z^d) )  $ or $ B_{2\delta_j} (y_j) \cap (\Delta_1 \cup \partial \eps ( D+\Z^d) ) = \emptyset $.

By the construction of $ \{ B_{\delta_j} (y_j) \}_{j=1}^{N_\eps} $, and the method we used in Proposition \ref{thm:micro-interior}, we only need to consider the case when $ y_j \in \partial\Omega $. In $ B_{2 \delta_j} (y_j) $, $ u_\eps $ solves
\begin{equation}
    -\diver(\mathbb C\ee(u_\eps))=0\ \text{in } B_{2\delta_j}(y_j)\cap U_1(x_0),
  \quad \mathbb{C}\ee( u_\eps)N =g \ \text{on }B_{2\delta_j}(y_j ) \cap \Delta_1,
\end{equation}
it then follows from the standard $ L^\infty  $ estimate of elliptic equation with Neumann boundary condition that:
\begin{equation}
     \begin{aligned}
         &\norm{\nabla u_\eps}_{L^\infty(B_{\delta_j}(y_j))}\\
  &\quad\le C \bigg\{\frac{1}{\delta_j}\inf_{c\in\R^d}\left(\avgint_{B_{2 \delta_j }(y_j)}\abs{u_\eps-c}^2\,dx\right)^{1/2} + \norm g_{L^\infty(\Delta_{2\delta_j }(y_j))}
   +\delta_j^\sigma[g]_{C^{0,\sigma}(\Delta_{2\delta_j}(y_j))} \bigg\}.
     \end{aligned}
\end{equation}

\eqref{eq:micro-boundary-estimate} then follows from summing up the $ L^\infty  $ estimate on $ B_{\delta_j} (y_j) $, and the  fact that
\begin{equation}
    \frac{1}{\delta_j} \le  C \eps^{-1} \le C \kappa.
\end{equation}
\end{proof}

\subsection{Large-scale interior Lipschitz estimate}


We first treat the interior problem. At scales larger than \(\eps\), the deviation from a rigid-affine first-order field is measured by an excess involving
\[
b+M(x-x_0)+\eps\chi_i^\alpha(x/\eps) M_i^\alpha.
\]
The local \(\Ahom\)-harmonic approximation, together with the regularity of the constant-coefficient system, gives a one-step decay of this excess. A coercivity property of the parametrization controls the affine matrices selected at successive scales, and iteration carries the estimate down to a scale comparable with \(\eps\). The remaining microscopic scale is handled by local interface estimates, the \(W^{1,\infty}\) bound for the cell correctors, and an affine nondegeneracy argument that recovers the skew-symmetric rotation inside each rigid inclusion.

For any $ u_\eps \in H^1(B_R(x_0);\R^d) $, such that $ \ee  (u_\eps ) = 0 $ in $ B_R(x_0) \cap \eps (D+\Z^d ) $, define:
\begin{equation}\label{eq:excess}
  H_{\eps}(u_\eps;x_0,r)
  :=\frac1r\inf_{\substack{b\in\R^d\\ M\in\R^{d\times d}}} \left(\avgint_{B_r(x_0)} \abs{u_\eps- \big(b+M(x-x_0)+ \eps \chi_i^\alpha(\frac{x}{\eps}) M_i^\alpha\big)}^2  \,dx\right)^{1/2}.
\end{equation}

\begin{lemma}\label{lem:parameter-coercivity}
There exist constants $\eta_0\in(0,1)$ and $C>0$ depending only on $ d,D,\lambda,\mu $, such that, for $0<\eps\le\eta_0r$,
\begin{equation}\label{eq:parameter-coercivity}
  \abs b+r\abs M
  \le C\left(\avgint_{B_r(x_0)} \abs{b+M(x-x_0)+ \eps\chi_i^\alpha(\frac{x}{\eps}) M_i^\alpha}^2\,dx\right)^{1/2}.
\end{equation}

Furthermore, assume for some $ u \in H^1(B_R;\R^d) $ with $ \ee(u) = 0 $ in $ \eps(D+\Z^d) $, there exist $ P_j = b_j + M_j(x-x_0 ) + \eps \chi_i^\alpha(\frac{x}{\eps}) (M_j)_i^\alpha$, and constants $ E_j>0 $, $ j=1,2 $, satisfying
\[
  \left(\avgint_{B_r(x_0)}\abs{u-P_j}^2\right)^{1/2}\le E_jr,
  \qquad j=1,2,
\]
then for $0<\eps\le\eta_0r$,
\begin{equation}\label{eq:adjacent-slope}
  \abs{b_1-b_2}+r\abs{M_1-M_2}
  \le C(E_1+E_2)r.
\end{equation}
\end{lemma}

\begin{proof}
   Denote $ P(x): = b+M(x-x_0)+ \eps \chi_i^\alpha(\frac{x}{\eps}) M_i^\alpha $ for short. It is clear that for $ h(x) : = b+ M(x-x_0) $, we have
    \begin{equation}
        \begin{aligned}
            |b|+ r|M|\le& C\left(\avgint_{B_r(x_0)} \abs{h}^2\,dx\right)^{1/2} \\
            \le & C \left(\avgint_{B_r(x_0)} \abs{P}^2\,dx\right)^{1/2} + C \eps\left(\avgint_{B_{r/\eps}(x_0/\eps)} \abs{\chi_i^\alpha M_i^\alpha}^2\,dx\right)^{1/2}, 
        \end{aligned}
    \end{equation}
and the last term of right hand side $ \left(\avgint_{B_{r/\eps}(x_0/\eps)} \abs{\chi_i^\alpha M_i^\alpha}^2\,dx\right)^{1/2} \le C |M| $, the first statement follows immediately from taking $ \eta_0 = \frac{1}{2C} $. The second statement follows just from triangle inequality.
\end{proof}

To prove the large-scale interior Lipschitz estimate, we need the following one-step improvement Lemma:

\begin{lemma}\label{thm:one-step}
Fix $ \nu \in (0,\beta) $, there exists a constant $ C=C(\lambda,\mu,d,D) $ such that for each $\theta\in(0,1/4]$, there is $\eps_1 = \eps_1(\theta)$ such that, for any $0<\eps\le\eps_1(\theta)r$ and $ u_\eps \in V_\eps(B_r(x_0)) $ satisfying
\begin{equation}
    \int_{B_r (x_0)} \mathbb{C} \ee(u_\eps) : \ee(\varphi) = 0, \quad \forall \varphi \in V_\eps(B_r(x_0))\cap H^1_0(B_r(x_0)),
\end{equation}
here $ V_\eps(B_r(x_0)) $ is defined by \eqref{def_VepsBall}, we have:
\begin{equation}\label{eq:one-step}
  H_{\eps}(u_\eps;x_0,\theta r)
  \le C \theta^\nu H_{\eps}(u_\eps ; x_0,r).
\end{equation}
Consequently, there exists $\theta_{0}\in(0,1/4]$ such that for each $ \eps\le\eps_1(\theta_0)r $,
\begin{equation}\label{eq:one-step-half}
  H_{\eps}(u_\eps; x_0,\theta_{0}r) \le\frac12 H_{\eps}(u_\eps ;x_0,r).
\end{equation}
\end{lemma}

\begin{proof}
Since \eqref{eq:one-step} is scale-invariant, we may assume $ r=1 $ without loss of generality. Assume also $ H_\eps (u_\eps ;x_0,1)>0 $, take $ P_0 = b_0 + M_0(x-x_0)+ \eps \chi_i^\alpha(\frac{x}{\eps}) (M_0)_i^\alpha $ such that
\begin{equation}
    \left(\avgint_{B_1(x_0)}\abs{u_\eps-P_0}^2\right)^{1/2}\le 2 H_\eps (u;x_0,1),
\end{equation}
and let $ v_\eps := (u_\eps(\cdot+x_0) -P_0) / (2H_\eps (u;x_0,1))$, then $ v_\eps \in H^1(B_1;\R^d) $, $ \ee(v_\eps) =0 $ in $ B_1\cap\eps(D+\Z^d) $, with $ \left(\avgint_{B_1}\abs{v_\eps}^2\right)^{1/2}\le 1 $, solves
\begin{equation}
    \mathrm{div} (\mathbb C \ee(v_\eps) ) = 0 \ \mathrm{in} \ B_1 \cap\eps(Y\setminus D + \Z^d).
\end{equation}
Here we used for each $ \phi \in V_\eps(B_1) \cap H^1_0(B_1) $, 
\begin{equation}
    \int_{B_1} \mathbb{C} \ee(M_0x+ \eps \chi_i^\alpha(\frac{x}{\eps}) (M_0)_i^\alpha) :\ee(\phi) = \int_{B_1} \widehat{\Sigma}_i^\alpha(x/\eps) (M_0)_i^\alpha :\ee(\phi) = 0.
\end{equation}

By Proposition \ref{thm:local-compactness}, take first $ \eps_1< \eps_0 $, then there exists an $\Ahom$-harmonic map $v$ in $B_{3/4}$ such that 
\begin{equation}
    \norm{v_\eps-v}_{L^2(B_{3/4})}\le C \eps^{\frac{1}{2}}, \quad \forall \eps\le \eps_1.
\end{equation}
By constant-coefficient interior regularity, there exists $h(x)=b+Mx$, such that for any fixed $\nu\in(0,1)$, 
\begin{equation}\label{eq:A-harmonic-affine}
  \left(\avgint_{B_\theta}\abs{v-h}^2\right)^{1/2}
  \le C\theta^{1+\nu},
  \qquad
  \abs b+\abs M\le C \|v\|_{L^2(B_{3/4})} \le C,
\end{equation}
for each $ \theta\in (0,1/2) $.
Let $ P (x) = b + Mx + \eps \chi_i^\alpha(\frac{x}{\eps}) M_i^\alpha $, we then derive
\begin{equation}
    \begin{aligned}
        & H_{\eps}(v_\eps; 0,\theta ) \\
        \le & \frac{1}{\theta} \left(\avgint_{B_\theta}\abs{v_\eps-P}^2\right)^{1/2} \\
        \le &\frac{1}{\theta} \left(\avgint_{B_\theta}\abs{v-h}^2\right)^{1/2} + \theta^{-(1+ \frac{d}{2})}\|v- v_\eps\|_{L^2(B_{3/4})} + \frac{1}{\theta} \left(\avgint_{B_\theta}\abs{\eps \chi_i^\alpha(x/\eps) M_i^\alpha }^2\right)^{1/2}\\
        \le & C ( \theta^\nu + \eps^{1/2} \theta^{-1-d/2} + \eps \theta^{-1} ).
    \end{aligned} 
\end{equation}

Take $ \eps_1 (\theta)\le \theta^{2(1+d/2+\nu)} $, we then derive \eqref{eq:one-step}.
\end{proof}



We are now in the position of proving the large-scale interior Lipschitz estimate:

Theorem \ref{thm_interior_lip_est} is a direct result of the following:

\begin{theorem}[Uniform interior Lipschitz estimate]\label{thm:interior-Lipschitz}
Assume $  \partial D$ is $ C^{2,\beta} $, there exists a constant $C= C(d,\mu,\lambda,D)$ such that, for any $\eps>0$, any ball $B_{2R}(x_0)$, and any $ u_\eps \in V_\eps(B_{2R}(x_0)) $ satisfying
\begin{equation}
    \int_{B_{2R} (x_0)} \mathbb{C} \ee(u_\eps) : \ee(\varphi) = 0, \quad \forall \varphi \in V_\eps(B_{2R}(x_0))\cap H_0^1(B_{2R} (x_0);\R^d),
\end{equation}
we have
\begin{equation}\label{eq:interior-Lipschitz}
  \norm{\nabla u_\eps}_{L^\infty(B_R(x_0))} \le\frac CR\inf_{c\in\R^d} \left(\avgint_{B_{2R}(x_0)}\abs{u_\eps-c}^2\,dx\right)^{1/2}.
\end{equation}
\end{theorem}

\begin{proof}
Set
\[
  h(u;x,\rho):=\inf_{c\in\R^d} \left(\avgint_{B_{\rho}(x)}\abs{u-c}^2\,dx\right)^{1/2},
\]
it is clear that $ H_\eps(u;x,\rho)\le\frac{1}{\rho} h(u;x,\rho) $.

Without loss of generality, assume $2d\eps<R $. Otherwise, use Proposition \ref{thm:micro-interior}. Fix $y\in B_R(x_0)$, $B_R(y)\subset B_{2R}(x_0)$. Take $ J \in \mathbb N_+$, such that $ \theta_0^{J+1} R \le \max \{\eta_0^{-1}, \eps_0^{-1}\}\eps < \theta_0^{J} R $, here $ \eta_0,\eps_0,\theta_0 $ are chosen as in Lemma \ref{lem:parameter-coercivity} and \ref{thm:one-step}. For each $ 0 \le j \le J$, let $ r_j = \theta_0^j R $, take $ P_j = b_j + M_j (x-y)+ \eps \chi_i^\alpha(\frac{x}{\eps}) (M_j)_i^\alpha$ such that
\begin{equation}\label{large_lip_est_eq1}
    \frac{1}{r_j} \left(\avgint_{B_{r_j}(y)} \abs{u_\eps-P_j}^2\right)^{1/2} \le  2 H_{\eps}(u_\eps;y,r_j ).
\end{equation}
since for each $ j\le J $, $\eps\le\eps_0 r_j$, by Lemma \ref{thm:one-step} we have:
\begin{equation}\label{large_lip_est_eq2}
    H_{\eps}(u_\eps;y,r_j)\le2^{-j}H_{\eps}(u_\eps;y,R).
\end{equation}

\emph{Step 1.} We now prove
\begin{equation}
  |M_J|\le C\big(\abs{M_0}+H_{\eps}(u_\eps;y,R)\big)\le C R^{-1} h(u_\eps;y,R).
\end{equation}
Notice that \eqref{large_lip_est_eq1} also implies that: for each $ j $,
 \begin{equation}
     \frac{1}{r_{j+1}} \left(\avgint_{B_{r_{j+1}}(y)} \abs{u_\eps -P_j}^2\right)^{1/2} \le  2 \theta_0^{-1-d/2} H_{\eps}(u_\eps;y,r_j ),
 \end{equation}
since for each $ 0 \le j\le J-1 $, $ \eps\le \eta_0 r_j $, Lemma \ref{lem:parameter-coercivity} then gives
\begin{equation}
    |b_j- b_{j+1}| + r_{j+1 }|M_j -M_{j+1}| \le C r_{j+1} \big(H_{\eps}(u_\eps;y,r_j )+H_{\eps}(u_\eps ;y,r_{j+1} )\big),
\end{equation}
recall \eqref{large_lip_est_eq2}, summing up gives
\begin{equation}
    |M_J| \le |M_0|+ C \sum_{j=0}^{J-1} 2^{-j}H_{\eps}(u_\eps ;y,R) \le|M_0|+ C H_\eps(u_\eps;y,R).
\end{equation}
It remains to estimate $ |M_0| $. Take $ c \in \R^d $ such that $ \left(\avgint_{B_{R}(y)}\abs{u_\eps-c}^2\,dx\right)^{1/2}  = h(u_\eps ;y,R) $, then again by Lemma \ref{lem:parameter-coercivity},
\begin{equation}
    \begin{aligned}
        R|M_0| &\le C\left(\avgint_{B_{R}(y)}\abs{P_0-c}^2\,dx\right)^{1/2}\\
        &\le 2 R CH_\eps(u_\eps;y,R) + Ch(u_\eps;y,R) \le Ch(u_\eps;y,R).
    \end{aligned}
\end{equation}
 
\emph{Step 2.} Let $w=u_\eps -P_J$. Then Proposition \ref{thm:micro-interior} and Lemma \ref{thm:one-step} give
\begin{equation}
    \begin{aligned}
        \|\nabla w\|_{L^\infty(B_{\frac{1}{2} \theta_0^J R } (y))} &\le \frac{C}{\theta_0^J R}\left(\avgint_{B_{\theta_0^J R}(y)} \abs w^2\right)^{1/2}\\
        &\le 2C  H_{\eps}(u_\eps; y,\theta_0^J R )\le C 2^{-J} C H_{\eps}(u_\eps; y, R ) \le \frac{C}{R} h(u_\eps;y,R),
    \end{aligned}
\end{equation}
combine with Lemma \ref{thm:corrector-W1infty} below, we then derive for each $ y\in B_{R}(x_0) $,
\begin{equation}
    |\nabla u_\eps (y)| \le  \|\nabla w\|_{L^\infty(B_{\frac{1}{2} \theta_0^J R } (y))} + |\nabla P_J| \le \frac{C}{R} h(u_\eps;y,R)\le \frac{C}{R} h(u_\eps;x_0,2R).
\end{equation}
\end{proof}

\begin{lemma}
    \label{thm:corrector-W1infty}
Assume $ \partial D \in C^{2,\beta} $, then there exists a constant $ C = C(d,D,\mu,\lambda) $, such that for each $ 1\le i,\alpha\le d $,
\begin{equation}\label{eq:corrector-W1infty}
  \norm{\nabla\chi_i^\alpha}_{L^\infty(\T^d)}\le C.
\end{equation}
\end{lemma}
\begin{proof}
Since $ \ee(\chi_i^\alpha + P_i^\alpha)  = 0 $ in $D$, there exists a rigid motion $ r_i^\alpha (y)= a_i^\alpha +B_i^\alpha y $, $ B_i^\alpha $ skew-symmetric, such that $ \chi_i^\alpha = -P_i^\alpha + r_i^\alpha $. By \eqref{eq:cell-bound},
\begin{equation}
    |B_i^\alpha| \le \|\nabla\chi_i^\alpha\|_{L^2(D)} +C |e^\alpha\otimes e^i| \le C,
\end{equation}
and
\begin{equation}
    |a_i^\alpha|\le |D|^{-1} \|r_i^\alpha\|_{L^2(D)}+d|B_i^\alpha|\le C|e^\alpha\otimes e^i| + \|\chi_i^\alpha\|_{L^2(D)} \le C.
\end{equation}

Since $ \chi_i^\alpha \in H^1(\T^d) $ solves
\begin{equation}
    \diver \mathbb{C}\ee(\chi_i^\alpha ) = 0 \ \mathrm{in} \ Y\setminus D, \quad \chi_i^\alpha|_{\partial D} =-P_i^\alpha + r_i^\alpha\in C^\infty(\partial D) ,
\end{equation}
it then follows by standard elliptic estimate with Dirichlet boundary condition that
\begin{equation}
    \|\nabla\chi_i^\alpha\|_{L^\infty(Y\setminus D)}\le C \|-P_i^\alpha+r_i^\alpha\|_{C^{1}} \le C.
\end{equation}
\end{proof}

\subsection{Large-scale boundary Lipschitz estimate}


We finally turn to the exterior traction boundary. The boundary excess measures both the deviation of the displacement from an affine field and the difference between the prescribed traction and the corresponding homogenized affine conormal. Boundary regularity for the constant-coefficient Neumann problem gives decay of this excess, while the local homogenized comparison transfers the decay to \(u_\eps\) with an error of order \((\eps/r)^{1/2}\). A discrete Campanato iteration controls the affine slopes and the displacement oscillation down to the microscopic scale. The microscopic boundary estimate closes the argument near \(\partial\Omega\), and the interior Lipschitz estimate controls the points lying outside this final boundary layer.

For $M\in\R^{d\times d}$, define:
\begin{equation}\label{eq:affine-conormal}
  \tau_M(x):=\Ahom M N(x),
  \qquad x\in\partial\Omega.
\end{equation}
and set
\begin{equation}\label{eq:boundary-flatness}
\begin{split}
  \Phi (u;x_0,r):= \inf_{\substack{b\in\R^d\\ M\in\R^{d\times d}}}
  \bigg\{&\frac1r \left(\avgint_{U_r}
       \abs{u-M(x-x_0)-b}^2\right)^{1/2}\\
  &+\norm{g-\tau_M}_{L^\infty(\Delta_r)}
   +r^\sigma[g-\tau_M]_{C^{0,\sigma}(\Delta_r)}\bigg\}.
\end{split}
\end{equation}

\begin{lemma}[Constant-coefficient Neumann flatness decay]\label{lem:hom-boundary-flatness}
There is $\theta_1\in ( 0,1/2) $, such that every solution $w$ of \eqref{eq:boundary-hom-problem} in $U_r$ satisfies
\begin{equation}\label{eq:hom-boundary-flatness}
  \Phi(w;x_0,\theta_1 r)\le\frac12\Phi(w; x_0,r).
\end{equation}
\end{lemma}

\begin{proof}
After translation and scaling, take $x_0=0$ and $r=1$. Take $ h_0(x) = M_0 x +b_0 $ such that
\begin{equation}
 \left(\avgint_{U_1}
       \abs{w-h_0}^2\right)^{1/2}+\norm{g-\tau_{M_0}}_{L^\infty(\Delta_1)}
   + [g-\tau_{M_0}]_{C^{0,\sigma}(\Delta_1)} \le 2 \Phi(w; 0,1).
\end{equation}
Let $ v:=w-h_0 $, then
\begin{equation}
    -\diver(\Ahom\ee(v))=0\quad\text{in }U_1,
  \qquad
  \Ahom\ee(v)N=g- \tau_{M_0}\quad\text{on }\Delta_1,
\end{equation}
standard boundary elliptic estimate gives: 
\begin{equation}
  \norm v_{C^{1,\sigma}(\overline{U_{1/2}})}
  \le C\Phi(w;0,1).
\end{equation}
Take $ M_1 = M_0+ \nabla v(0) $, $ b_1= b_0 +v(0)$, then we have:
\begin{equation}\label{eq:hom-flatness-oscillation}
  \theta^{-1}
  \left(\avgint_{U_{\theta}}
    \abs{w-M_1x-b_1}^2\right)^{1/2}
  \le C\theta^{\sigma}\Phi(w;0,1).
\end{equation}

For the boundary terms, define:
\begin{equation}
    \tau_v(x):=\Ahom\big(\ee(v)(x)-\ee(v)(0)\big),
\end{equation}
then $ g- \tau_{M_1} = \tau_v N $ on $ \Delta_1 $. Hence,
\begin{equation}
    \norm{g-\tau_{M_1}}_{L^\infty(\Delta_{\theta})} \le C\sup_{x\in \Delta_\theta}|\ee(v)(x)-\ee(v)(0)|
    \le C\theta^{\sigma}\Phi(w;0,1),
\end{equation}
and for $x,y\in\Delta_{\theta}$,
\begin{equation}
    \begin{aligned}
        &\abs{\tau_v(x)N(x)-\tau_v(y)N(y)}\\
        \le & \abs{\tau_v(x)-\tau_v(y)}
     +\abs{\tau_v(y)}\abs{N(x)-N(y)}\\
     \le & C |\ee(v)(x)-\ee(v)(y)|+ C\|e(v)-e(v)(0)\|_{L^\infty(\Delta_\theta)} \abs{N(x)-N(y)}\\
     \le & C(|x-y|^\sigma + \theta^\sigma |x-y|^{\beta})\Phi(w;0,1).
    \end{aligned}
\end{equation}
here for the third inequality, we used $N\in C^{0,\beta}$. Therefore taking $ \sigma\le \beta $ gives
\begin{equation}\label{eq:hom-flatness-traction-holder}
  \theta^{\sigma} [g-\tau_{M_1}]_{C^{0,\sigma}(\Delta_{\theta})}
  \le C\big(\theta^{\sigma}+\theta^{\sigma+\beta}\big)\Phi(w;0,1).
\end{equation}
Combine all above, we then derive:
\begin{equation}
    \begin{aligned}
        \theta^{-1}
  \left(\avgint_{U_{\theta}}
    \abs{w-M_1x-b_1}^2\right)^{1/2}+\norm{g-\tau_{M_1}}_{L^\infty(\Delta_{\theta})} + \theta^{\sigma} [g-\tau_{M_1}]_{C^{0,\sigma}(\Delta_{\theta})}\qquad \\
    \le C \theta^\sigma \Phi(w; 0,1),
    \end{aligned}
\end{equation}
\eqref{eq:hom-boundary-flatness} then follows from choosing $ \theta_1\in(0,1/2) $ such that $ C\theta_1^\sigma<\frac{1}{2} $.
\end{proof}

Define
\begin{equation}\label{eq:Phi-boundary}
  \Phi^{\partial}(u_\eps;x_0,r):=\frac1r\inf_{c\in\R^d}
       \left(\avgint_{U_{2r}}\abs{u_\eps-c}^2\right)^{1/2}
       +\norm g_{L^\infty(\Delta_{2r})}.
\end{equation}

\begin{lemma}[One-step boundary inequality]\label{lem:boundary-one-step}
Let $\theta_1$ be as in Lemma \ref{lem:hom-boundary-flatness}, there exists $C<\infty$ such that, for any $ \eps \le c_0 r $,
\begin{equation}\label{eq:boundary-one-step}
  \Phi(u_\eps; x_0,\theta_1 r)
  \le\frac12\Phi(u_\eps;x_0,r)
     +C\left(\frac\eps r\right)^{1/2}\Phi^{\partial}(u_\eps;x_0,r).
\end{equation}
\end{lemma}

\begin{proof}
Let $u_0$ be the solution from \cref{prop:boundary-comparison}. Then by Lemma \ref{lem:hom-boundary-flatness},
\begin{align*}
  \Phi (u_\eps; x_0,\theta_1 r)
  &\le\Phi(u_0;x_0,\theta_1 r)
      +\frac C r\left(\avgint_{U_r}\abs{u_\eps-u_0}^2\right)^{1/2}\\
  &\le\frac12 \Phi(u_0;x_0,r)
      + \frac {C} r\left(\avgint_{U_{r}}\abs{u_\eps-u_0}^2\right)^{1/2}\\
  &\le\frac12\Phi(u_\eps ; x_0,r)
      +2\frac C r\left(\avgint_{U_{r}}\abs{u_\eps-u_0}^2\right)^{1/2},
\end{align*}
\eqref{eq:boundary-one-step} then follows from  \eqref{eq:boundary-comparison}.
\end{proof}

\begin{lemma}[Discrete boundary iteration]\label{lem:boundary-iteration}
There exists $\kappa_1 \ge 4d$, such that the following holds: For any $ R> \kappa_1\eps $, take $ J \in \mathbb{N} $ such that $ \theta_1^{J+1} R< \kappa_1 \eps \le \theta_1^J R $, let $r_j=\theta_{1}^jR$. For each $j$, assume $ Q_j(x) : = M_j(x-x_0)+b_j $ satisfies
\begin{equation}
     \frac1{r_j} \left(\avgint_{U_{r_j}} \abs{u_\eps-Q_j}^2\right)^{1/2}
  +\norm{g-\tau_{M_j}}_{L^\infty(\Delta_{r_j})}
   +r_j^\sigma[g-\tau_{M_j}]_{C^{0,\sigma}(\Delta_{r_j})} \le 2 \Phi (u_\eps;x_0,r_j),
\end{equation}
Then
\begin{equation}\label{eq:boundary-iteration-bound}
  \max_{0\le j\le J}
  \set{\Phi (u_\eps; x_0,r_j)+\abs{M_j}}
  \le C\phi(u_\eps;x_0,3R/2),
\end{equation}
where
\begin{equation}\label{eq:Lambda-R}
  \phi(u_\eps ;x_0,R):=\frac1R\inf_{c\in\R^d}
       \left(\avgint_{U_{R}}\abs{u_\eps-c}^2\right)^{1/2}
       +\norm g_{L^\infty(\Delta_{R})}
       +R^\sigma[g]_{C^{0,\sigma}(\Delta_{R})}.
\end{equation}
Moreover, for every $r\in[\kappa_{1}\eps,R]$,
\begin{equation}\label{eq:boundary-oscillation-all-scales}
  \frac1r\inf_{c\in\R^d}
       \left(\avgint_{U_r}\abs{u_\eps-c}^2\right)^{1/2}
  \le C\phi(u_\eps ;x_0,3R/2 ).
\end{equation}
\end{lemma}

\begin{proof}
Choose first $\kappa_{1}\ge c_{0}^{-1}$ with $ c_0 $ defined as in Lemma \ref{lem:boundary-one-step}. 
If $ J=0 $, then by $ N\in C^{0,\beta}(\Delta_R) , \sigma<\beta $,
\begin{equation}
    \begin{aligned}
        \Phi(u_\eps;x_0,R) & \le \phi (u_\eps;x_0,3R/2),
    \end{aligned}
\end{equation}
hence, by triangle inequality,
\begin{equation}\label{one_step_lem_eq4}
    \begin{aligned}
        |M_0| & \le C  \frac1R\inf_{c\in\R^d} \left(\avgint_{U_R}\abs{ Q_0-c}^2\right)^{1/2} \\
        &\le C \Phi(u_\eps;x_0,R)+C \phi(u_\eps;x_0,R) \le C \phi (u_\eps;x_0,3R/2).
    \end{aligned}
\end{equation}

We now assume $ J\ge 1 $. Since for each $ j $, $ r_j \ge c_0^{-1}\eps $, Lemma \ref{lem:boundary-one-step} gives
\begin{equation}\label{one_step_lem_eq2}
    \Phi (u_\eps; x_0,r_{j+1})
  \le\frac12\Phi(u_\eps;x_0,r_j)
     +C\left(\frac\eps r_j\right)^{1/2}\Phi^{\partial}(u_\eps; x_0,r_j).
\end{equation}
For $ j\ge1 $, by $ \theta_1\le1/2 $, $ U_{2r_j} \subset U_{r_{j-1}} $, we have:
\begin{equation}\label{one_step_lem_eq3}
    \begin{aligned}
        \Phi^{\partial}(u_\eps;x_0,r_j) \le &  \frac1{r_j}\inf_{c\in\R^d}
       \left(\avgint_{U_{2r_j}}\abs{u_\eps-c-M_{j-1} x}^2\right)^{1/2}
       +\norm g_{L^\infty(\Delta_{2r_j})}+ C|M_{j-1}| \\
       \le & 2 \Phi(u_\eps;x_0,r_{j-1}) + C|M_{j-1}|,
    \end{aligned}
\end{equation}
and for $ j=0 $, by definition, $ \Phi^{\partial}(u_\eps;x_0,R) \le \phi(u_\eps; x_0,R)$. 

Since it is clear that the second part of Lemma \ref{lem:parameter-coercivity} still holds for $ P $ replaced by affine maps, we have:
\begin{equation}
    r_{j+1} |M_{j+1} -M_j | +|b_{j+1} -b_j | \le C r_{j+1} \big(\Phi (u_\eps;x_0,r_{j+1})+\Phi( u_\eps; x_0,r_{j}) \big),
\end{equation}
Summing up gives 
\begin{equation}\label{one_step_lem_eq1}
    \begin{aligned}
        \max_{0\le j\le J}|M_j| \le |M_0| + C \sum_{j=0}^{J}  \Phi(u_\eps; x_0,r_{j}),
    \end{aligned}
\end{equation}
using \eqref{one_step_lem_eq2} and \eqref{one_step_lem_eq3} gives
\begin{equation}
    \begin{aligned}
        \sum_{j=0}^{J}  \Phi(u_\eps; x_0,r_{j})
        \le &\phi(u_\eps;x_0,3R/2) + \frac{1}{2}\sum_{j=0}^{J}  \Phi(u_\eps;x_0,r_{j} ) \\
        & \quad+ C \sum_{j=1}^{J} \left(\frac\eps r_j\right)^{1/2}\big(\Phi(u_\eps;x_0,r_{j-1}) + |M_{j-1}|\big).
    \end{aligned}
\end{equation}
Denote $ S_J = \sum_{j=0}^{J}  \Phi( u_\eps; x_0,r_{j})  $, then
\begin{equation}
    \begin{aligned}
        S_J & \le 2 \phi(u_\eps;x_0,3R/2)+ C \kappa_1^{-1/2} S_J + C \eps^{1/2}R^{-1/2} (|M_0|+S_J) \sum_{j=1}^{J} \theta_1^{-j/2}  \\
        & \le 2 \phi(u_\eps;x_0,3R/2)+ C \kappa_1^{-1/2} (S_J + |M_0|),
    \end{aligned}
\end{equation}
here we used \eqref{one_step_lem_eq1}, the assumption $ \eps/r_j <\kappa_1^{-1} $ for each $ j $ and $ R^{-1/2}\sum_{j=1}^{J-1} \theta_1^{-j/2} \le C (\theta_1^{J} R)^{-1/2} \le \eps^{-1/2} \kappa_1^{-1/2} $. 

Now take $ \kappa_1^{1/2} \ge 2C $, we then derive by \eqref{one_step_lem_eq4}:
\begin{equation}
    S_J \le 4 \phi(u_\eps;x_0,3R/2)+ |M_0| \le C \phi(u_\eps;x_0,3R/2).
\end{equation}
In particular, $  \Phi(u_\eps;x_0,r_j) \le S_J \le C \phi(u_\eps;x_0,3R/2)$.

For \eqref{eq:boundary-oscillation-all-scales}, take $ 0\le j\le J $ such that $ \theta_1r_j \le r<r_j $, then
\begin{equation}
    \begin{aligned}
        \frac1r\inf_{c\in\R^d} \left(\avgint_{U_r}\abs{u_\eps-c}^2\right)^{1/2} & \le C \frac1{r_j} \inf_{c\in\R^d} \left(\avgint_{U_{r_j}}\abs{u_\eps-c}^2\right)^{1/2}\\
        & \le C \big( \Phi(u_\eps;x_0,r_j) + |M_j|\big) \le  C \phi(u_\eps;x_0,3R/2).
    \end{aligned}
\end{equation}
\end{proof}

\begin{proof}[Proof of Theorem~\ref{thm:local-boundary-Lipschitz}]
Without loss of generality, assume $ \kappa_1\eps<R $. Otherwise, use Proposition \ref{thm:micro-boundary}.
Take $ J\in \mathbb N_+ $ such that $ \theta_1^{J+1} R \le \kappa_1\eps < \theta_1^{J} R $. Fix $ y \in \Delta_{\frac{R}{2}}(x_0) $, for each $ 0\le j \le J$, let $ r_j =\theta_1^j R $ take $ Q_j (x)= b_j + M_j(x-y) $, such that
\begin{equation}\label{thm_boundary_lip_eq2}
    \begin{aligned}
        \frac1{r_j} \left(\avgint_{U_{r_j} (y)}
       \abs{u_\eps -Q_j}^2\right)^{1/2}+\norm{g-\tau_{M_j}}_{L^\infty(\Delta_{r_j })}  +r_j^\sigma[g-\tau_{M_j}]_{C^{0,\sigma}(\Delta_{r_j})} \qquad \\
       \le 2 \Phi (u_\eps; y,r_j),
    \end{aligned}
\end{equation}
then by Lemma \ref{lem:boundary-iteration}, we have
\begin{equation}\label{thm_boundary_lip_eq1}
  \begin{aligned}
      \max_{0\le j\le J}\Phi(u_\eps;y,r_j )+\abs{M_j} + \sup_{r\in[\kappa_1 \eps, R]} \frac1r\inf_{c\in\R^d}
       \left(\avgint_{U_r}\abs{u_\eps-c}^2\right)^{1/2}\qquad\\
       \le C\phi(u_\eps;y,3R/2) \le C \phi(u_\eps;x_0,2R).
  \end{aligned}
\end{equation}

Now taking $ \kappa= \kappa_1\theta_1^{-1} $ in Proposition \ref{thm:micro-boundary} gives:
\begin{equation}
    \begin{aligned}
        \norm{\nabla u_\eps}_{L^\infty(U_{\frac{1}{2} r_J}(y))}
 & \le C \phi (u_\eps; y,r_J ) \\
 & \le C(2 \Phi (u_\eps; y,r_J) + |M_J| ) \le C\phi(u_\eps;x_0,2R),
    \end{aligned}
\end{equation}
here we used $ N\in C^{0,\beta}(\Delta_R) , \sigma<\beta$. Hence, we have:
\begin{equation}
    \norm{\nabla u_\eps}_{L^\infty(U_{\frac{1}{2}r_J}(y))} \le C\phi(u_\eps;y,R)+|M_J| \le C\phi(u_\eps;x_0, 2R).
\end{equation}

We now show that $|\nabla u_\eps(x)|$ is controlled by the right-hand side of \eqref{eq:local-boundary-Lipschitz} for all $x\in U_{R/2}(x_0)$ lying in the boundary layer, namely, for all points covered by \[ \bigcup_{y\in\Delta_{R/2}(x_0)} U_{r_J/2}(y). \] 

For the remaining point $ y $ in $ U_{R/2} (x_0) $, denote $ d_y : = \mathrm{dist} (y,\Delta_{R/2}(x_0)) = |y- y_0|$ for some $ y_0\in \Delta_{R/2}(x_0) $, then $ d_y \ge \frac{1}{2} r_J > \frac{1}{2}\kappa_1\eps >2d\eps$, by the choice of $ \kappa_1 $ in Lemma \ref{lem:boundary-iteration}. Hence, $ \mathrm{dist}  (B(y,\frac{1}{2}d_y) ,\Delta_{R/2}(x_0))> d \eps $, it follows $ B(y,\frac{1}{2}d_y) \cap D_\eps = B(y,\frac{1}{2}d_y) \cap \eps(D+\Z^d) $. 

We can now use Theorem \ref{thm:interior-Lipschitz} on $ u_\eps $ to give:
\begin{equation}
    \begin{aligned}
        \norm{\nabla u_\eps}_{L^\infty(B_{d_y/4} (y))} \le& \frac C{d_y}\inf_{c\in\R^d} \left(\avgint_{B_{d_y/2} (y)}\abs{u_\eps-c}^2\,dx\right)^{1/2}\\
        \le & \frac C{d_y}\inf_{c\in\R^d} \left(\avgint_{U_{3d_y/2} (y_0)}\abs{u_\eps-c}^2\,dx\right)^{1/2}.
    \end{aligned}
\end{equation}
If $ \frac{3}{2} d_y \in[\kappa_1\eps, R] $, use \eqref{thm_boundary_lip_eq1} to give:
\begin{equation}
    \begin{aligned}
        \norm{\nabla u_\eps}_{L^\infty(B_{d_y/4} (y))} \le C \phi(u_\eps;x_0,2R) .
    \end{aligned}
\end{equation}
If $ \frac{3}{2}d_y>R $, we simply use:
\begin{equation}
     \frac 1{d_y}\inf_{c\in\R^d} \left(\avgint_{U_{3d_y/2} (y_0)}\abs{u_\eps-c}^2\,dx\right)^{1/2} \le  C \phi(u_\eps;x_0,2R) .
\end{equation}
 If $ \frac{3}{2}d_y < \kappa_1 \eps $, then $ \kappa_1\eps /2 < d_y <2\kappa_1\eps/3 $, $ d_y $ is comparable to $ \kappa_1\eps $, we have
\begin{equation}
    \begin{aligned}
        \norm{\nabla u_\eps}_{L^\infty(B_{d_y/4} (y))} \le \frac C{\kappa_1\eps }\inf_{c\in\R^d} \left(\avgint_{U_{\kappa_1 \eps} (y_0)}\abs{u_\eps-c}^2\,dx\right)^{1/2} \le C \phi(u_\eps;x_0,2R).
    \end{aligned}
\end{equation}
\end{proof}

\section*{Acknowledgements}

The authors thank Xin Fu, Wenjia Jing, and Beichen Wang for their guidance
and suggestions. This research benefited from substantial assistance from
AIM (AI Mathematician)~\cite{liu2026aimathematicianfullyautomated}, an
AI-assisted mathematical research system. Through repeated interactions
with the authors, AIM contributed to the progressive clarification of the
research questions, the development of comparison constructions that
preserve the rigidity of the inclusions, and the exploration of arguments
connecting quantitative homogenization, the control of microscopic
rotations, and boundary regularity. The authors carefully verified and
revised these contributions and take full responsibility for the
mathematical statements and proofs presented in this paper.

\bibliographystyle{plainnat}
\bibliography{reference2}

\end{document}